\documentclass[12pt,aps,tightenlines]{revtex4}
\pdfoutput=1

\usepackage[letterpaper]{geometry}
\usepackage{amsmath}
\usepackage{amsfonts}
\usepackage{bm}
\usepackage{graphicx}
\usepackage{subfigure}
\usepackage{color}

\graphicspath{{figs_pdf/}}

\newcommand{\Sobo}{W}

\newcommand{\Small}{\delta}
\newcommand{\Stresstensor}{T}
\newcommand{\eye}{i}

\newcommand{\pairing}[2]{\left(\!\left(#1\,,\,#2\right)\!\right)}

\begin{document}
%+Title
\title{Nonlinear dynamics of phase separation in ultra-thin films}
\author{Lennon \'O N\'araigh}
\affiliation{Department of Mathematics, Imperial College
  London, SW7 2AZ, United Kingdom}
\author{Jean-Luc Thiffeault}
\email{jeanluc@mailaps.org}
\affiliation{Department of Mathematics, University of Wisconsin, Madison,
WI 53706, USA}
\date{7 May 2008}
%-Title
%
\begin{abstract}
  We present a long-wavelength approximation to the Navier--Stokes
  Cahn--Hilliard equations to describe phase separation in thin films.
  The equations we derive underscore the coupled behaviour of
  free-surface variations and phase separation.  Since we are
  interested in the long-time behaviour of the phase-separating fluid,
  we restrict our attention to films that do not rupture.  To do this,
  we introduce a regularising Van der Waals potential.  We analyse the
  resulting fourth-order equations by constructing a solution as the
  limit of a Galerkin approximation, and obtain existence and
  regularity results.  In our analysis, we find a nonzero lower bound
  for the height of the film, which precludes the possibility of
  rupture.  The lower bound depends on the parameters of the problem,
  and we compare this dependence with numerical simulations.  We find
  that while the theoretical lower bound is crucial to the
  construction of a smooth, unique solution to the PDEs, it is
  not sufficiently sharp to represent accurately the parametric
  dependence of the observed dips in free-surface height.
\end{abstract}

\maketitle
\section{Introduction}
\label{sec:intro} 

\noindent Below a certain critical temperature, a well-mixed binary
fluid spontaneously separates into its component parts, forming
domains of pure liquid.  This process can be characterised by the
Cahn--Hilliard equation, and numerous studies describe the physics and
mathematics of phase separation.  In this paper we study phase
separation in a thin layer, in which the the varying free surface and
concentration fields are coupled through a pair of nonlinear evolution
equations.

Cahn and Hilliard introduced their eponymous equation
in~\cite{CH_orig} to model phase separation in a binary alloy.  Since
then, the model has been used in diverse applications: to describe
polymeric fluids~\cite{Aarts2005}, fluids with interfacial
tension~\cite{LowenTrus}, and self-segregating populations in
biology~\cite{Murray1981}.
Analysis of the Cahn--Hilliard (CH) equation was given by Elliott and
Zheng in~\cite{Elliott_Zheng}, where they prove the existence,
uniqueness, and regularity of solutions: given sufficiently smooth
initial data, solutions to the CH equation are bounded for almost all
time in the Sobolev space $\Sobo^{4,2}$.  Several authors have
developed generalisations of the CH equation: a variable-mobility
model was introduced by Elliott and Garcke~\cite{Elliott_varmob},
while nonlocal effects were considered by Gajweski and
Zacharias~\cite{Gajewski_nonlocal}.  These additional features do not
qualitatively change the phase separation, and we therefore turn to
one mechanism that does: the coupling of a flow field to the
Cahn--Hilliard equation~\cite{Bray_advphys}.
In this case, the Cahn--Hilliard concentration equation is modified by an
advection term, and the flow field is either prescribed or evolves according
to some fluid equation.  Ding and co-workers~\cite{Spelt2007} provide
a derivation of coupled Navier--Stokes Cahn--Hilliard (NSCH) equations in
which the velocity advects the phase-separating concentration field, while
concentration gradients modify the velocity through an additional stress
term in the momentum equation.  Such models have formed the basis of numerical
studies of binary fluids~\cite{Berti2005}, while other studies without
this feedback term highlight different regimes of phase separation under
flow~\cite{chaos_Berthier, ONaraigh2006, Lacasta1995}.  In this paper, the
NSCH equations form the starting point for our asymptotic analysis.

As in other applications involving the Navier--Stokes equations, the
complexity of the problem is reduced when the fluid is spread thinly
on a substrate, and the upper vertical boundary forms a free
surface~\cite{Oron1997}.  Then, provided vertical gradients are small
compared to lateral gradients, a long-wavelength approximation is
possible, in which the full equations with a moving boundary at the
free surface are reduced to a single equation for the free-surface
height.  In the present case, the reduction yields two equations: one
for the free surface, and one for the Cahn--Hilliard concentration.
The resulting thin-film Stokes Cahn--Hilliard equations have already
been introduced by the authors in~\cite{ONaraigh2007}, although the
focus there was on control of phase separation and numerical
simulations in three dimensions.  Here we confine ourselves to the
two-dimensional case: we derive the thin-film equations from first
principles, present analysis of the resulting equations, and highlight
the impossibility of film rupture.
 
Along with the simplification of the problem that thin-film theory
provides, there are many practical reasons for studying phase
separation in thin layers.  Thin polymer films are used in the
fabrication of semiconductor devices, for which detailed knowledge of
film morphology is required~\cite{Karim2002}.  Other industrial
applications of polymer films include paints and coatings, which are
typically mixtures of polymers.  One potential application of the
thin-film Cahn--Hilliard theory is in
self-assembly~\cite{Putkaradze2005,Xia2004, Krausch1994}.  Here
molecules (usually residing in a thin layer) respond to an
energy-minimisation requirement by spontaneously forming large-scale
structures.  Equations of Cahn--Hilliard type have been proposed to
explain the qualitative features of
self-assembly~\cite{Putkaradze2005,Holm2005}, and knowledge of
variations in the film height could enhance these models.  Indeed
in~\cite{ONaraigh2007} the authors use the present thin-film
Cahn--Hilliard model in three dimensions to control phase separation,
a useful tool in applications where it is necessary for the molecules
in the film to form a given structure.

The analysis of thin-film equations was given great impetus by Bernis
and Friedman in~\cite{Friedman1990}.  They focus 
on the basic thin-film equation,
\begin{equation}
\frac{\partial h}{\partial t}=-\frac{\partial}{\partial x}\left(h^n\frac{\partial^3{h}}{\partial{x^3}}\right),
\label{eq:basic_h}
\end{equation}
with no-flux boundary conditions on a line segment, and smooth
nonnegative initial conditions.  For $n=1$ this equation describes a
thin bridge between two masses of fluid in a Hele--Shaw cell, for
$n<3$ it is used in slip models as $h\rightarrow0$~\cite{Myers1998},
while for $n=3$ it gives the evolution of the free surface of a thin
film experiencing capillary forces~\cite{Oron1997}.
Using a decaying free-energy functional, they prove the existence of
nonnegative solutions to Eq.~\eqref{eq:basic_h} for $n\geq1$, while
for $n\geq 4$, and for strictly positive initial conditions, the solution
is unique, strictly positive, and is almost
always bounded in the $\Sobo^{3,2}$ norm.  This paper has inspired
other work on the subject~\cite{Bertozzi1996,
  Bertozzi1998,Laugesen2002}, in which the effect of a Van der Waals
term on Eq.~\eqref{eq:basic_h} is investigated.  These works provide
results concerning regularity, long-time behaviour, and film rupture
in the presence of an attractive Van der Waals force.
More relevant to the present work is the paper by Wieland and
Garcke~\cite{Wieland2006}, in which a pair of partial differential
equations describes the coupled evolution of free-surface variations
and surfactant concentration.  The authors derive the relevant
equations using the long-wavelength theory, obtain a decaying energy
functional, and prove results concerning the existence and
nonnegativity of solutions.  We shall take a similar approach in this
paper.

When the binary fluid forms a thin film on a substrate, we shall show
in Sec.~\ref{sec:model} that a long-wave approximation simplifies the
Navier--Stokes Cahn--Hilliard equations, which reduce to a pair of
coupled evolution equations for the free surface and concentration.
If $h(x,t)$ is the scaled free-surface height, and
$c(x,t)$ is the binary fluid concentration, then the
dimensionless equations take the form
\begin{subequations}
\begin{equation}
\frac{\partial h}{\partial t}+\frac{\partial J}{\partial x}=0,\qquad
\frac{\partial}{\partial t}\left(h c\right)+\frac{\partial}{\partial x}\left(Jc\right)=\frac{\partial}{\partial{x}}\left(h\frac{\partial\mu}{\partial{x}}\right),
\end{equation}
where
\begin{equation}
J = \tfrac{1}{2} h^2 \frac{\partial\sigma}{\partial{x}}
- \tfrac{1}{3} h^3 \bigg\{
           \frac{\partial}{\partial{x}}
           \left(-\frac{1}{C}\frac{\partial^2{h}}{\partial{x}^2}
             + \phi\right) + \frac{r}{h} \frac{\partial}{\partial{x}}
           \left[h{\left(\frac{\partial{c}}{\partial{x}}
             \right)}^2\right]\bigg\},
\end{equation}
\begin{equation}
\mu=c^3-c-\frac{C_{\mathrm{n}}^2}{h}\frac{\partial}{\partial{x}}\left(h\frac{\partial{c}}{\partial{x}}\right).
\end{equation}%
\label{eq:model_intro}%
\end{subequations}%
Here $C$ is the capillary number, $r$ measures the strength of
coupling between the concentration and free-surface variations
(backreaction), and $C_{\mathrm{n}}$ is the scaled interfacial
thickness.  Additionally, $\sigma$ is the dimensionless,
spatially-varying surface tension, and $\phi$ is the body-force
potential acting on the film.  In this paper we take
$\phi=-|A|h^{-3}$, the repulsive Van der Waals
potential~\cite{Book_Parsegian2006}.  This choice stabilises the film
and prevents rupture.  Although rupture is in itself an important
feature in thin-film equations~\cite{Oron1997, Bertozzi1996,
  Bertozzi1998}, in this paper we are interested in late-time phase
separation and it is therefore undesirable.

We present the asymptotic analysis that converts the NSCH equations
into Eq.~\eqref{eq:model_intro} and prove that the model equations
possess smooth solutions that are bounded for almost all time in the
$\Sobo^{4,2}$ norm.  The principal tool in this analysis is the
construction of a free-energy functional for
Eq.~\eqref{eq:model_intro} that is a decaying function of time.  We
prove that $h(x,t)>0$ for all time, which is the no-rupture condition.

The paper is organised as follows.  In Sec.~\ref{sec:model} we discuss
the Navier--Stokes Cahn--Hilliard equation and the scaling laws that
facilitate the passage to the long-wavelength equations, and we derive
Eq.~\eqref{eq:model_intro}.  In Sec.~\ref{sec:existence} we analyse
these equations by constructing a decaying free-energy functional.  We
prove the existence of solutions to Eq.~\eqref{eq:model_intro} and
provide regularity results.  We obtain a condition on the minimum
value of the free-surface height, and show that this is never zero.  Using
numerical studies, we discuss the
dependence of the minimum free-surface height on the problem
parameters in Sec.~\ref{sec:height_dip}, and compare with the analytic
results.  Finally, in Sec.~\ref{sec:conclusions} we present our
conclusions.

\section{The Model Equations}
\label{sec:model}
\noindent 
In this section we introduce the two-dimensional Navier--Stokes
Cahn--Hilliard (NSCH) equation set.
We discuss the assumptions underlying the long-wavelength
approximation.
We enumerate the scaling rules necessary to obtain the simplified
equations.  
Finally, we arrive at a set of equations that describe phase separation in
a thin film subject to arbitrary body forces.

The full NSCH equations describe the coupled effects of phase separation
and flow in a binary fluid.  If $\bm{v}$ is the fluid velocity and $c$ is
the concentration of the mixture, where $c=\pm1$ indicates total segregation,
then these fields evolve as
\begin{subequations}
\begin{gather}
        \frac{\partial \bm{v}}{\partial t}+\bm{v}\cdot\nabla\bm{v}=\nabla\cdot\Stresstensor-\frac{1}{\rho}\nabla\phi,\\
        \frac{\partial c}{\partial t}+\bm{v}\cdot\nabla c=D \nabla^2\left(c^3-c-\gamma\nabla^2 c\right),\\
        \nabla\cdot\bm{v}=0,
\end{gather}%
\label{eq:NSCH}%
\end{subequations}%
where
\begin{equation}
        \Stresstensor_{ij} =-\frac{p}{\rho}\delta_{ij}+\nu\left(\frac{\partial
        v_i}{\partial
        x_j}+\frac{\partial
        v_j}{\partial x_i}\right)-\beta\gamma\frac{\partial c}{\partial x_i}\frac{\partial
        c}{\partial x_j}
\label{eq:NSCH_tensor}%
\end{equation}%
is the stress tensor, $p$ is the fluid pressure, $\phi$ is the body
potential and $\rho$ is the constant density.  The constant $\nu$ is
the kinematic viscosity, $\nu=\eta/\rho$, where $\eta$ is the dynamic
viscosity.  Additionally, $\beta$ is a constant with units of
$[\mathrm{Energy}][\mathrm{Mass}]^{-1}$, $\sqrt{\gamma}$ is a constant
that gives the typical width of interdomain transitions, and $D$ a
diffusion coefficient with dimensions
$[\mathrm{Length}]^2[\mathrm{Time}]^{-1}$.

If the system has a free surface in the vertical or $z$-direction and has
infinite or periodic boundary conditions (BCs) in the lateral or $x$-direction,
then the vertical BCs we impose are 
\begin{subequations}
\begin{equation}
u = w = c_z = c_{zzz}\text{ on }z=0,
\end{equation}
while on the free surface $z=h(x,t)$ they are
\begin{equation}
\hat{n}_i \hat{n}_j \Stresstensor_{ij} = -\sigma\kappa,\qquad \hat{n}_i \hat{t}_j
\Stresstensor_{ij} = -\frac{\partial\sigma}{\partial
s},
\label{eq:BC_stress}
\end{equation}
\begin{equation}
w=\frac{\partial h}{\partial t}+u\frac{\partial h}{\partial x},
\label{eq:BC_height}
\end{equation}
\begin{equation}
\hat{n}_i \partial_i c = 0, \qquad \hat{n}_i \partial_i \nabla^2 c = 0,
\label{eq:BC_c}
\end{equation}%
\label{eq:BCs}%
\end{subequations}%
where $\hat{\bm{n}} =
(-\partial_x{h}\,,\,1)/{[{1+(\partial_x{h})^2}]^{1/2}}$ is the unit
normal to the surface, $\hat{\bm{t}}$, is the unit vector tangent to
the surface, $s$ is the surface coordinate, $\sigma$ is the surface
tension, and $\kappa$ is the mean curvature,
\begin{equation*}
\kappa=\nabla\cdot\hat{\bm{n}}
%=\left(\partial_x,\partial_z\right)\cdot\left(\frac{-\partial_x{h}}{\sqrt{1+\left(\partial_x{h}\right)^2}},\frac{1}{\sqrt{1+\left(\partial_x{h}\right)^2}}\right)
=\frac{\partial_{xx}h}{\left[1+\left(\partial_x{h}\right)^2\right]^{\frac{3}{2}}}.
\end{equation*}
This choice of BCs guarantees the conservation of the total mass and volume,
\begin{equation}
\text{Mass} = \int_{\text{Dom}(t)} c(\bm{x},t)d^2x,\qquad
\text{Volume} = \int_{\text{Dom}(t)}d^2x.
\label{eq:integral_quantities}
\end{equation}
Here 
$\text{Dom}(t)%=\{\left(x,y\right)|x\in\left[0,L\right],z\in\left[0,h(x,t)\right]\}
$
represents the time-dependent domain of integration, owing to the
variability of the free surface height.  Note that in view of the
concentration BC~\eqref{eq:BC_c}, the stress BC~\eqref{eq:BC_stress}
and does not contain $c(\bm{x},t)$ or its derivatives.

These equations simplify considerably if the fluid forms a thin layer
of mean thickness $h_0$, for then the scale of lateral variations
$\lambda$ is large compared with the scale of vertical variations
$h_0$.  Specifically, the parameter $\Small = h_0/\lambda$ is small,
and after nondimensionalisation of Eq.~\eqref{eq:NSCH} we expand its
solution in terms of this parameter, keeping only the lowest-order
terms.  For a review of this method and its applications,
see~\cite{Oron1997}.  For simplicity, we shall work in two dimensions,
but the generalisation to three dimensions is easily
effected~\cite{ONaraigh2007}.

In terms of the small parameter $\Small$, the equations
nondimensionalise as follows.  The diffusion time scale is $t_0 =
\lambda^2 / D = h_0^2 / \left(\Small^2 D\right)$ and we choose this to
be the unit of time.  Then the unit of horizontal velocity is $u_0 =
\lambda/ t_0 = \Small D / h_0$ so that $u = \left(\Small D /
  h_0\right)U$, where variables in upper case denote dimensionless
quantities.  Similarly the vertical velocity is $w = \left(\Small^2 D
  / h_0\right)W$.  For the equations of motion to be half-Poiseuille
at $O\left(1\right)$ (in the absence of the backreaction) we choose $p
= \left(\eta D / h_0^2\right)P$ and $\phi = \left(\eta D /
  h_0^2\right)\Phi$.  Using these scaling rules, the dimensionless
momentum equations are
\begin{multline}
\Small Re \left(\frac{\partial U}{\partial T} + U\frac{\partial
U}{\partial X}+ W\frac{\partial U}{\partial Z}\right) 
= -\frac{\partial}{\partial X}\left(P+\Phi\right)+\Small^2\frac{\partial^2
U}{\partial X^2}+\frac{\partial^2 U}{\partial Z^2} \\
-\tfrac{1}{2}\frac{\beta\gamma}{\nu D}\frac{\partial}{\partial
X}\bigg[\Small^2\left(\frac{\partial c}{\partial X}\right)^2
+\left(\frac{\partial c}{\partial Z}\right)^2\bigg]-\frac{\beta\gamma}{\nu
D}\frac{\partial c}{\partial X}\bigg[\Small^2\frac{\partial^2 c}{\partial X^2}+\frac{\partial^2 c}{\partial Z^2}\bigg],
\label{eq:thin_film1}
\end{multline}
\begin{multline}
\Small^3 Re \left(\frac{\partial W}{\partial T} + U\frac{\partial
W}{\partial X} + W\frac{\partial W}{\partial Z}\right) 
= -\frac{\partial}{\partial Z}\left(P+\Phi\right)+\Small^4\frac{\partial^2
W}{\partial X^2}+\Small^2\frac{\partial^2 W}{\partial Z^2}  \\
-\tfrac{1}{2}\frac{\beta\gamma}{\nu D}\frac{\partial}{\partial
Z}\bigg[\Small^2\left(\frac{\partial c}{\partial X}\right)^2+\left(\frac{\partial
c}{\partial Z}\right)^2\bigg]
-\frac{\beta\gamma}{\nu D}\frac{\partial c}{\partial Z}\bigg[\Small^2\frac{\partial^2
c}{\partial X^2}+\frac{\partial^2 c}{\partial Z^2}\bigg],
\label{eq:thin_film2}
\end{multline}
where
\begin{equation}
Re = \frac{\Small D}{\nu} = \frac{\left(\Small D / h_0\right)h_0}{\nu} =
O\left(1\right).
\end{equation}
The choice of ordering for the Reynolds number $Re$ allows us to recover
half-Poiseuille flow at $O\left(1\right)$.  We delay choosing the ordering
of the dimensionless group $\beta\gamma/D\nu$ until we have examined the
concentration equation, which in nondimensional form is
\begin{multline}
\Small^2 \left(\frac{\partial c}{\partial T} + U\frac{\partial c}{\partial
X} + W\frac{\partial c}{\partial Z}\right)\\
= \Small^2 \frac{\partial^2}{\partial X^2}\left(c^3-c\right)+\frac{\partial^2}{\partial Z^2}\left(c^3-c\right)
-\Small^4 C_{\mathrm{n}}^2\frac{\partial^4c}{\partial X^4}-C_{\mathrm{n}}^2\frac{\partial^4
c}{\partial Z^4}
-2\Small^2 C_{\mathrm{n}}^2 \frac{\partial^2}{\partial X^2}\frac{\partial
c}{\partial Z^2},
\label{eq:conc0}
\end{multline}
where $C_{\mathrm{n}}^2=\gamma/h_0^2$.  By switching off the
backreaction in the momentum equations (corresponding to
$\beta\gamma/D\nu\rightarrow\infty$), we find the trivial solution to
the momentum equations, $U = W = \partial_X\left(P+\Phi\right)
= \partial_Z\left(P+\Phi\right) = 0$, $H = 1$.  The concentration
boundary conditions are then $c_Z = c_{ZZZ} = 0$ on $Z = 0,1$ which
forces $c_Z\equiv0$ so that the Cahn--Hilliard equation is simply
\[
\frac{\partial c}{\partial T} = \frac{\partial^2}{\partial X^2}\left(c^3-c\right)-\Small^2 C_{\mathrm{n}}^2\frac{\partial^4c}{\partial
X^4}.
\]
To make the lubrication approximation consistent, we take 
\begin{equation}
\Small C_{\mathrm{n}}= \tilde{C}_{\mathrm{n}}= \Small\sqrt{\gamma} /
h_0 = O\left(1\right).
\label{eq:gamma_h}
\end{equation}
We now carry out a long-wavelength approximation to
Eq.~\eqref{eq:conc0}, writing $U=U_0+O\left(\Small\right)$,
$W=W_0+O\left(\Small\right)$, $c=c_0+\Small c_1+\Small^2 c_2+\ldots$.
We examine the boundary conditions on $c(\bm{x},t)$ first.
They are $\hat{\bm{n}}\cdot\nabla c =
\hat{\bm{n}}\cdot\nabla\nabla^2c=0$ on $Z=0,H$; on $Z=0$ these
conditions are simply $\partial_Z c = \partial_{ZZZ} c=0$, while on
$Z=H$ the surface derivatives are determined by the relations
\[
\hat{\bm{n}}\cdot\nabla\ \propto\ -\Small^2 H_X\partial_X+\partial_Z,
\]
\[
\hat{\bm{n}}\cdot\nabla\nabla^2\ \propto\ -\Small^4
H_X\partial_{XXX}-\Small^2H_X\partial_X\partial_{ZZ}+\Small^2\partial_{XX}\partial_Z+\partial_{ZZZ}.
\]
Thus, the BCs on $c_0$ are simply $\partial_Z c_0 = \partial_{ZZZ} c_0=0$
on $Z=0,H$, which forces $c_0=c_0\left(X,T\right)$.  Similarly, we
find  $c_1=c_1\left(X,T\right)$, and
\[
\frac{\partial c_2}{\partial Z}=Z\frac{H_X}{H}\frac{\partial
c_0}{\partial X},\qquad
\frac{\partial^2 c_2}{\partial Z^2}=\frac{H_X}{H}\frac{\partial
c_0}{\partial X}, \qquad \text{for any }Z\in \left[0,H\right].
\]
In the same manner, we derive the results $\partial_{ZZZZ}c_2=\partial_{ZZZZ}c_3=0$.
Using these facts, Eq.~\eqref{eq:conc0} becomes
\begin{multline*}
\frac{\partial c_0}{\partial T} + U_0\frac{\partial c_0}{\partial
X} =
\\
\frac{\partial^2}{\partial X^2}\left(c_0^3-c_0\right)
-\tilde{C}_{\mathrm{n}}^2\frac{\partial^4c}{\partial X^4}
+\left(3c_0^2-1\right)\frac{H_X}{H}\frac{\partial c_0}{\partial X}
-2\tilde{C}_{\mathrm{n}}^2\frac{\partial^2}{\partial X^2}\frac{H_X}{H}\frac{\partial
c_0}{\partial X}
-{\tilde{C}_{\mathrm{n}}^2}\frac{\partial^4 c_4}{\partial Z^4}.
\end{multline*}
We now integrate this equation from $Z=0$ to $H$ and use the boundary conditions
\[
\frac{\partial^3 c_4}{\partial Z^3}=0\quad\text{on}\quad Z=0,
\]
\begin{equation*}
\frac{\partial^3 c_4}{\partial Z^3}=
H_X\frac{\partial^3 c_0}{\partial X^3}+
H_X\frac{\partial}{\partial X}\left(\frac{H_X}{H}\frac{\partial c_0}{\partial
X}\right)
-H\frac{\partial^2}{\partial X^2}\left(\frac{H_X}{H}\frac{\partial c_0}{\partial
X}\right)\quad\text{on}\quad Z=H.
\end{equation*}
After rearrangement, the concentration equation becomes
\begin{multline*}
H\frac{\partial c_0}{\partial T}+H\langle U_0\rangle\frac{\partial c_0}{\partial
X} = 
\\
H\frac{\partial^2}{\partial X^2}\left[c_0^3-c_0-\tilde{C}_{\mathrm{n}}^2\frac{\partial^2
c_0}{\partial X^2}-\tilde{C}_{\mathrm{n}}^2\frac{H_X}{H}\frac{\partial c_0}{\partial
X}\right]
+\frac{\partial H}{\partial X}\frac{\partial}{\partial X}\left[c_0^3-c_0-\tilde{C}_{\mathrm{n}}^2\frac{\partial^2
c_0}{\partial X^2}-\tilde{C}_{\mathrm{n}}^2\frac{H_X}{H}\frac{\partial c_0}{\partial
X}\right],
\end{multline*}
where
\[
\langle U_0\rangle = \frac{1}{H}\int_0^H U_0\left(X,Z,T\right)dZ
\]
is the vertically-averaged velocity.  Introducing
\[
\mu = c_0^3-c_0-\frac{\tilde{C}_{\mathrm{n}}^2}{H}\frac{\partial}{\partial X}\left(H\frac{\partial c_0}{\partial
X}\right),
\]
the thin-film Cahn--Hilliard equation becomes
\begin{equation}
\frac{\partial c_0}{\partial T}+\langle U_0\rangle\frac{\partial c_0}{\partial
X} = \frac{1}{H}\frac{\partial}{\partial X}\left(H\frac{\partial\mu}{\partial
X}\right).
\label{eq:conc_eqn_ornament}
\end{equation}

We are now able to perform the long-wavelength approximation to Eqs.~\eqref{eq:thin_film1}
and~\eqref{eq:thin_film2}.  At lowest order, Eq.~\eqref{eq:thin_film2} is~$\partial_ Z\left(P+\Phi\right)=0$,
since~$c_0=c_0(X,T)$, and hence
\[
P+\Phi = P_{\mathrm{surf}}+\Phi_{\mathrm{surf}}\equiv P\left(X,h(x,t),T\right)+\Phi\left(X,h(x,t),T\right).
\]
By introducing the backreaction strength
\begin{equation}
r=\frac{\Small^2\beta\gamma}{D\nu}=O\left(1\right),
\end{equation}
equation~\eqref{eq:thin_film1} becomes
\[
\frac{\partial^2 U_0}{\partial Z^2}=\frac{\partial}{\partial X}\left(P_{\mathrm{surf}}+\Phi_{\mathrm{surf}}\right)+r\frac{\partial}{\partial
X}\left(\frac{\partial c_0}{\partial X}\right)^2+r\frac{\partial c_0}{\partial
X}\frac{\partial^2 c_2}{\partial Z^2}.
\]
Using $\partial_{ZZ}c_2=\left(H_X/H\right)\left(\partial c_0/\partial
X\right)$ this becomes
\begin{equation}
\frac{\partial^2 U_0}{\partial Z^2}=\frac{\partial}{\partial X}\left(P_{\mathrm{surf}}+\Phi_{\mathrm{surf}}\right)+\frac{r}{H}\frac{\partial}{\partial
X}\left[H\left(\frac{\partial c_0}{\partial X}\right)^2\right].
\label{eq:d2U}
\end{equation}
At lowest order, the BC~\eqref{eq:BC_stress} reduces to
\[
\frac{\partial U_0}{\partial
  Z}=\frac{\partial\Sigma}{\partial{X}}\quad\text{on}\quad Z=H,
\]
which combined with Eq.~\eqref{eq:d2U} yields the relation
\[
\frac{\partial U_0}{\partial Z}=\frac{\partial\Sigma}{\partial X}+\left(Z-H\right)\bigg\{
\frac{\partial}{\partial X}\left(P_{\mathrm{surf}}+\Phi_{\mathrm{surf}}\right)+\frac{r}{H}\frac{\partial}{\partial
X}\left[H\left(\frac{\partial c_0}{\partial X}\right)^2\right]
\bigg\}.
\]
Here $\Sigma$ is the dimensionless, spatially-varying surface tension.  Making
use of the BC $U_0=0$ on $Z=0$ and integrating again, we obtain the result
\begin{equation}
U_0\left(X,Z,T\right)=Z\frac{\partial\Sigma}{\partial X}+\left(\tfrac{1}{2}Z^2-HZ\right)\bigg\{
\frac{\partial}{\partial X}\left(P_{\mathrm{surf}}+\Phi_{\mathrm{surf}}\right)+\frac{r}{H}\frac{\partial}{\partial
X}\left[H\left(\frac{\partial c_0}{\partial X}\right)^2\right]
\bigg\}.
\label{eq:U_0}
\end{equation}
The vertically-averaged velocity is therefore
\[
\langle U_0\rangle = \tfrac{1}{2}H\frac{\partial\Sigma}{\partial X}-\tfrac{1}{3}H^2\bigg\{
\frac{\partial}{\partial X}\left(-\frac{1}{C}\frac{\partial^2 H}{\partial X^2}+\Phi_{\mathrm{surf}}\right)+\frac{r}{H}\frac{\partial}{\partial
X}\left[H\left(\frac{\partial c_0}{\partial X}\right)^2\right]
\bigg\},
\]
where we used the standard Laplace--Young free-surface boundary
condition to eliminate the pressure, and
\begin{equation}
C = \frac{\nu\rho D}{h_0\sigma_0\Small^2}=O\left(1\right).
\label{eq:C}
\end{equation}
Finally, by integrating the continuity equation in the $Z$-direction, we
obtain, in a standard manner, an equation for free-surface variations,
\begin{equation}
\frac{\partial H}{\partial X}+\frac{\partial}{\partial X}\left(H\langle U_0\rangle\right)=0.
\label{eq:h_eqn_ornament}
\end{equation}

Let us assemble our results, restoring the lower-case fonts and omitting
ornamentation over the constants.  The height equation~\eqref{eq:h_eqn_ornament}
becomes
\begin{subequations}
%\begin{gather}
\begin{equation}
\frac{\partial h}{\partial t}+\frac{\partial J}{\partial x}=0,
\end{equation}
while the concentration equation~\eqref{eq:conc_eqn_ornament} becomes
\begin{equation}
\frac{\partial}{\partial t}\left(c h\right)+\frac{\partial}{\partial x}\left(Jc\right)=\frac{\partial}{\partial{x}}\left(h\frac{\partial\mu}{\partial{x}}\right),
\end{equation}
where
\begin{equation}
J=\tfrac{1}{2}h^2\frac{\partial\sigma}{\partial{x}}-\tfrac{1}{3}h^3\bigg\{\frac{\partial}{\partial{x}}\left(-\frac{1}{C}\frac{\partial^2{h}}{\partial{x}^2}
+\phi\right)+\frac{r}{h}\frac{\partial}{\partial{x}}\left[h\left(\frac{\partial{c}}{\partial{x}}\right)^2\right]\bigg\},
\end{equation}
and
\begin{equation}
\mu=c^3-c-C_{\mathrm{n}}^2\frac{1}{h}\frac{\partial}{\partial{x}}\left(h\frac{\partial{c}}{\partial{x}}\right),
\end{equation}%
\label{eq:model}%
\end{subequations}%
and where we have the nondimensional constants
\begin{equation}
r=\frac{\Small^2\beta\gamma}{D\nu},\qquad C_{\mathrm{n}}=\frac{\Small\sqrt{\gamma}}{h_0},\qquad
C=\frac{\nu\rho D}{h_0\sigma_0\Small^2}.
\end{equation}
These are the thin-film NSCH equations.  The integral quantities defined in Eq.~\eqref{eq:integral_quantities} are manifestly conserved,
while the free surface and concentration are coupled.

We note that the relation
$C_{\mathrm{n}}=\delta\sqrt{\gamma}/h_0=O\left(1\right)$ is the
condition that the mean thickness of the film be much smaller than the
transition layer thickness.  In experiments involving the smallest
film thicknesses attainable ($10^{-8}$ m)~\cite{Sung1996}, this
condition is automatically satisfied.  The condition is also realised
in ordinary thin films when external effects such as the air-fluid and
fluid-substrate interactions do not prefer one binary fluid component
or another.  In this case, the vertical extent of the domains becomes
comparable to the film thickness at late times, the thin film behaves
in a quasi two-dimensional way, and the model equations are
applicable.

The choice of potential $\phi$ determines the behaviour of solutions.
If interactions between the fluid and the substrate and air interfaces
are important, the potential should take account of the Van der Waals
forces present.  A simple model potential is thus
\[
\phi = Ah^{-n},
\]
where $A$ is the dimensionless Hamakar constant and typically
$n=3$~\cite{Oron1997}.  Here $A$ can be positive or negative, with
positivity indicating a net attraction between the fluid and the
substrate and negativity indicating a net repulsion.  This choice of
potential can also have a regularising effect, preventing a
singularity or rupture from occurring in Eq.~\eqref{eq:model}.
  
For $\phi=-|A|/h^3$ (repulsive Van der Waals interactions), the system of
equations~\eqref{eq:model} has a Lyapunov functional $F=F_1+F_2$, where
\[
F_1=\int{dx}\, \left[\frac{1}{2C}\left(\frac{\partial{h}}{\partial{x}}\right)^2+\frac{|A|}{2h^2}\right],\qquad
F_2 = \frac{r}{C_{\mathrm{n}}^2}\int{dx}\,  h\left[\tfrac{1}{4}\left(c^2-1\right)^2+\frac{C_{\mathrm{n}}^2}{2}\left(\frac{\partial{c}}{\partial{x}}\right)^2\right].
\]  
By differentiating these expressions with respect to time, we obtain
the relation
\begin{multline}
\dot{F}_1+\dot{F}_2
\\
=-\tfrac{1}{3} \int{dx}\,  h^3\bigg\{\frac{\partial}{\partial x}\left(\frac{1}{C}\frac{\partial^2h}{\partial
x^2}+\frac{|A|}{h^3}\right)-\frac{r}{h}\frac{\partial}{\partial x}\left[h\left(\frac{\partial
c}{\partial x}\right)^2\right]  \bigg\}^2
-\int{dx}\,  h\left(\frac{\partial\mu}{\partial x}\right)^2,
\label{eq:fe_decay}
\end{multline}
which is nonpositive for nonnegative $h$.  This fact is the key to
the analytic results of the next section. 

\section{Existence of solutions to the model equations}
\label{sec:existence}

\noindent In this section we prove that solutions to the model equations
do indeed exist.  We set $C=\tfrac{1}{3}, r=|A|=1$ in Eqs.~\eqref{eq:model}
and focus on the resulting equation set
\begin{subequations}
\begin{equation}
\frac{\partial h}{\partial t}=
-\frac{\partial}{\partial x}\left[f(h)\frac{\partial^{3}h}{\partial{x^3}}\right]+
\frac{\partial}{\partial x}\left[\frac{1}{g(h)}\frac{\partial h}{\partial{x}}\right]+
\frac{\partial}{\partial x}\bigg\{\frac{f(h)}{g(h)}\frac{\partial}{\partial{x}}\left[g(h)\left(\frac{\partial{c}}{\partial{x}}\right)^2\right]\bigg\},
\end{equation}
%
%
%\vskip -0.2in
\begin{multline}
\frac{\partial}{\partial t}\left(c g(h)\right)=
-\frac{\partial}{\partial x}\left[cf(h)\frac{\partial^{3}h}{\partial{x^3}}\right]+
\frac{\partial}{\partial x}\left[\frac{c}{g(h)}\frac{\partial h}{\partial{x}}\right]+
\frac{\partial}{\partial x}\bigg\{c\frac{f(h)}{g(h)}\frac{\partial}{\partial{x}}\left[g(h)\left(\frac{\partial{c}}{\partial{x}}\right)^2\right]\bigg\}
\\
+\frac{\partial}{\partial x}\bigg\{g(h)\frac{\partial}{\partial
x}\left[c^3-c-\frac{1}{g(h)}\frac{\partial}{\partial
x}\left(g(h)\frac{\partial c}{\partial x}\right)\right]\bigg\},
\label{eq:model_proof_c}%
\end{multline}%
\label{eq:model_proof}%
\end{subequations}%
where
\[
f(h)=h^3,\qquad g(h)=h.
\]
The equations are defined on a periodic domain $\Omega=[0,L]$, while
the initial conditions are
\begin{equation}
h(x,0)=h_0(x)>0,\qquad c(x,0)=c_0(x),\qquad
h_0(x),\,c_0(x)\in \Sobo^{2,2}(\Omega),
\label{eq:initial_data}
\end{equation}
We shall prove that the solutions to this equation pair exist in the strong
sense; however, we shall need the definition of weak solutions:
\begin{quote}
A pair $\left(h,c\right)$ is a \emph{weak solution} 
of Eq.~\eqref{eq:model_proof} if the following integral relations hold:
\begin{multline*}
\int_0^{T_0}dt\int_\Omega dx\,\varphi_t h=
\\
\int_0^{T_0}dt\int_\Omega dx\,\varphi_x\bigg\{-f(h)\frac{\partial^{3}h}{\partial{x^3}}+
\frac{1}{g(h)}\frac{\partial h}{\partial{x}}+
\frac{f(h)}{g(h)}\frac{\partial}{\partial{x}}\left[g(h)\frac{\partial}{\partial{x}}\left(\frac{\partial{c}}{\partial{x}}\right)^2\right]\bigg\},
\end{multline*}
and
\begin{multline}
\int_0^{T_0}dt\int_\Omega dx\,\psi_t cg(h)=
\\
\int_0^{T_0}dt\int_\Omega dx\,\psi_x\bigg\{-cf(h)\frac{\partial^{3}h}{\partial{x^3}}+
\frac{c}{g(h)}\frac{\partial h}{\partial{x}}+
c\frac{f(h)}{g(h)}\frac{\partial}{\partial{x}}\left[g(h)\frac{\partial}{\partial{x}}\left(\frac{\partial{c}}{\partial{x}}\right)^2\right]\bigg\}
\\
+\int_0^{T_0}dt\int_\Omega dx\,\psi_x \bigg\{g(h)\frac{\partial}{\partial
x}\left[c^3-c-\frac{1}{g(h)}\frac{\partial}{\partial{x}}\left(g(h)\frac{\partial{c}}{\partial{x}}\right)\right]\bigg\},
\label{eq:weak_sln}
\end{multline}
where $T_0>0$ is any time, and $\varphi\left(x,t\right)$ and $\psi\left(x,t\right)$
are arbitrary differentiable test functions that are periodic on $\Omega$
and vanish at $t=0$ and $t=T_0$.  
\end{quote}
In a series of steps in Secs.~\ref{sec:regularization}--\ref{sec:uniqueness},
we prove this result:
\begin{quote}
  \emph{Given the initial data in Eq.~\eqref{eq:initial_data},
    Eqs.~\eqref{eq:model_proof} possess a strong solution endowed with
    the following regularity properties:
\[
\left(h,c\right)\in L^{\infty}\left(0,T_0;\Sobo^{2,2}(\Omega)\right)
\cap L^2\left(0,T_0;\Sobo^{4,2}(\Omega)\right)
\cap C^{\frac{3}{2},\frac{1}{8}}\left(\Omega\times\left[0,T_0\right]\right).
\]
}
\end{quote}
The outline of the proof is as follows: In
Sec.~\ref{sec:regularization} we introduce a regularised version of
Eqs.~\eqref{eq:model_proof}, whose solution we approximate by a
Galerkin sum in Sec.~\ref{sec:galerkin}.  In
Sec.~\ref{sec:a_priori_bds}, we obtain \emph{a priori} bounds on
various norms of the approximate solution.  Crucially, we show that
the free-surface height is always positive.  This enables us to
continue the approximate solution in the time interval
$\left[0,T_0\right]$.  In Secs.~\ref{sec:equicontinuity}
and~\ref{sec:convergence} we show that the Galerkin sum converges to a
solution of the unapproximated equations, in the appropriate limit.
Finally, in Secs.~\ref{sec:regularity} and~\ref{sec:uniqueness} we
discuss the regularity and uniqueness properties of the solution.

\subsection{Regularisation of the problem}
\label{sec:regularization}

\noindent We introduce regularised functions
$f_\varepsilon(s)$ and $g_\varepsilon(s)$ such
that $\lim_{\varepsilon\rightarrow0}f_\varepsilon(s)=
f(s)$, and
$\lim_{\varepsilon\rightarrow0,s\geq0}g_\varepsilon(s)={g}(s)$.
For now we do not specify $f_\varepsilon(s)$, although we
mention that a suitable choice of $f_\varepsilon(s)$ will
cure the degeneracy of the fourth-order term in the height equation.
On the other hand, we require that~$g_\varepsilon(s)$ have the properties:
(i) $g_\varepsilon(s)=s+\varepsilon$, for $s\geq0$;
(ii) $g_\varepsilon(s)>0$, for $s<0$;
(iii) $\lim_{s\rightarrow-\infty}g_\varepsilon(s)=\tfrac{1}{2}\varepsilon$;
and (iv) $g_\varepsilon(s)$ is at least $C^3$.

From Eqs.~\eqref{eq:model_proof}, the regularised PDEs we study are
\begin{subequations}
\begin{align}
%\begin{split}
h_t&=-J_{\varepsilon,x}\,,
\\
\left(cg_\varepsilon(h)\right)_t&=-\left(c J_\varepsilon\right)_x-\left(g_\varepsilon(h)\mu_{\varepsilon,x}\right)_x,
\label{eq:reg_pde_c}%
%\end{split}%
\end{align}%
\label{eq:reg_pde}%
\end{subequations}%
where
\[
\mu_\varepsilon=c^3-c-\frac{1}{g_\varepsilon(h)}\left(g_\varepsilon(h)c_x\right)_x
\]
and
\[
J_\varepsilon=f_\varepsilon(h)h_{xxx}-\frac{1}{g_\varepsilon(h)}h_x-\frac{f_\varepsilon(h)}{g_\varepsilon(h)}\left(g_\varepsilon(h)c_{x}^2\right)_x\,.
\]
Equation~\eqref{eq:reg_pde_c} can also be written as 
\begin{equation}
c_t=-\frac{1}{g_\varepsilon(h)}J_\varepsilon c_x-\frac{1}{g_\varepsilon(h)}\left(g_\varepsilon(h)\mu_{\varepsilon,x}\right)_x-\frac{c}{g_\varepsilon(h)}\left[J_{\varepsilon,x}+g_\varepsilon'(h)h_t\right];
\label{eq:conc_no_flux}
\end{equation}
this form of the concentration equation will be useful in
Sec.~\ref{sec:uniqueness}.

\subsection{The Galerkin approximation}
\label{sec:galerkin}

\noindent We choose a complete orthonormal basis on the interval
$\Omega$, with periodic boundary conditions.  Let us denote the basis
by $\{\phi_i(x)\}_{i\in\mathbb{N}_0}$.  We consider the
finite vector space $\text{Span}\{\phi_0,\ldots,\phi_n\}$.  For
convenience, let us take the $\phi_i(x)$'s to be the
eigenfunctions of the Laplacian on $\left[0,L\right]$ with periodic
boundary conditions, and corresponding eigenvalues $-\lambda_i^2$.
Let $\phi_0$ be the constant eigenfunction.  We construct
approximate solutions to the PDEs~\eqref{eq:reg_pde} as finite sums,
\[
h_n(x,t)=\sum_{i=0}^n\eta_{n,i}(t)\phi_i(x),
\qquad c_n(x,t)=\sum_{i=0}^n\gamma_{n,i}(t)\phi_i(x).
\]
If the (smooth) initial data are given as
\[
h(x,0)=h_0(x)=\sum_{i=0}^\infty\eta_i^0\phi_i(x)>0,\qquad
c(x,0)=c_0(x)=\sum_{i=0}^\infty\gamma_i^0\phi_i(x),
\]
then the initial data for the Galerkin approximation are
\[
h_n(x,0)=h_n^0(x)=\sum_{i=0}^n\eta_{i}^0\phi_i(x),\qquad
c_n(x,0)=c_n^0(x)=\sum_{i=0}^n\gamma_{i}^0\phi_i(x),
\]
and the initial data of the Galerkin approximation converge strongly in the
$L^2(\Omega)$ norm to the initial data of the unapproximated problem.
Thus, there is a $n_0\in\mathbb{N}$ such that $h_n^0(x)>0$, everywhere
in $\Omega$, for all $n>n_0$.  Henceforth we work with Galerkin approximations
with $n>n_0$. 

Substitution of $h_n=\sum_{i=0}^n\eta_{n,i}\phi_i$  into a weak form of the
$h$-equation yields
\begin{equation}
\frac{d}{dt}\pairing{h_n}{\phi_j}=\pairing{J_{\varepsilon,n}}{\phi_{j,x}},
\label{eq:weak_h}
\end{equation}
where
\[
J_\varepsilon\left(h_n,c_n\right)=f_\varepsilon(h_n)h_{n,xxx}-\frac{1}{g_\varepsilon(h_n)}h_{n,x}-\frac{f_\varepsilon(h_n)}{g_\varepsilon(h_n)}\left(g_\varepsilon\left(h_{n}\right)c_{n,x}^2\right)_x,
\]
is the flux for the regularised $h$-equation, and $\pairing{\varphi(x)}{\psi(x)}$
is the pairing $\int_\Omega\varphi\psi\,{dx}$.
We recast Eq.~\eqref{eq:weak_h} as
\[
\frac{d\eta_{n,j}}{dt}=\pairing{J_{\varepsilon,n}}{\phi_{j,x}}=\Phi_{n,j}\left(\eta_n,\gamma_n\right),
\]
where the function $\Phi_n\left(\eta_n,\gamma_n\right)$ depends on
$\eta_n=\left(\eta_{n,0},\ldots,\eta_{n,n}\right)$ and
$\gamma_n=\left(\gamma_{n,0},\ldots,\gamma_{n,n}\right)$, and is
locally Lipschitz in its variables.  This Lipschitz property arises
from the fact that the regularised flux, evaluated at the Galerkin
approximation, is a composition of Lipschitz continuous functions, and
therefore, is itself Lipschitz continuous.

Similarly, substitution of $c_n=\sum_{i=0}^n\gamma_{n,i}\phi_i$ into the
weak form of the $c$-equation~\eqref{eq:reg_pde_c} yields
\begin{equation}
\frac{d}{dt}\pairing{g (h_n)c_n}{\phi_j}=\pairing{K_{\varepsilon,n}}{\phi_{j,x}},
\label{eq:weak_c}
\end{equation}
where
\begin{align*}
  K_\varepsilon\left(h_n,c_n\right) &=
  c_nf_\varepsilon(h_n)h_{n,xxx} -
  \frac{c_n}{g_\varepsilon(h_n)}h_x-c_{n}
  \frac{f_\varepsilon(h_n)} {g_\varepsilon(h_n)}
  \left(g_\varepsilon(h_n)
    c_{n,x}^2\right)_x - g (h_n)\mu_{\varepsilon,n,x}\,,\\
  &=
  c_nJ_{\varepsilon,n}-g_\varepsilon(h_n)\mu_{\varepsilon,n,x}\,,
\end{align*}
is the flux for the regularised
$c$-equation~\eqref{eq:reg_pde_c}. 
Rearranging gives
\[
\pairing{g (h_n)c_{n,t}}{\phi_j}=\pairing{K_{\varepsilon,n}}{\phi_{j,x}}-\pairing{g' (h_n)c_{n}h_{n,t}}{\phi_j},
\]
and the left-hand side can be recast in matrix form as
$\sum_{\eye=0}^nM_{\eye{j}}\dot{\gamma}_{n,\eye}$,
where
\[
M_{\eye j}=\int_\Omega{dx}\,  g_\varepsilon\left(\sum_\ell \eta_{n,\ell}\phi_\ell\right)\phi_\eye\phi_j,
\] 
which is manifestly symmetric.  It is positive
definite because given a vector $\left(\xi_0,\ldots,\xi_n\right)$, we
have the relation
\begin{equation*}
\sum_{\eye,j}\xi_\eye M_{\eye j}\xi_j=\int_\Omega{dx}\,  g_\varepsilon\left(\sum_\ell\eta_{n,\ell}\phi_\ell\right)\sum_{\eye,j}\left(\phi_\eye\xi_\eye\right)\left(\phi_j\xi_j\right)
>0,\quad\text{for }\left(\xi_0,\ldots,\xi_n\right)\neq0,
\end{equation*}
which follows from the positivity of the regularised function $g_\varepsilon(s)$.
 We therefore have the following equation for $\gamma_{n,j}(t)$,
\begin{equation}
\frac{d\gamma_{n,j}}{dt}=\sum_{\eye=0}^n M_{\eye j}^{-1}\Big[\pairing{K_{\varepsilon,n}}{\phi_{\eye,x}}-\pairing{g_\varepsilon' (h_n)c_{n}h_{n,t}}{\phi_\eye}
\Big].
\label{eq:gamma_dot}
\end{equation}
Inspecting the expression for $M_{\eye j}$, $g (h_n)$, and
$K_{\varepsilon,n}$, we see that $\eta_n$ and $\gamma_n$ appear in a
Lipschitz-continuous way in the expression $\sum_{\eye=0}^n M_{\eye
  j}^{-1}\pairing{K_{\varepsilon,n}}{\phi_{\eye,x}}$, while owing to
the imposed smoothness of $g_\varepsilon(s)$, the variables
$\eta_n$, $\gamma_n$ and
$\dot{\eta}_n=\left(\dot{\eta}_{n,0},\ldots,\dot{\eta}_{n,n}\right)$
appear in a Lipschitz-continuous way in the quantity $\sum_{\eye=0}^n
M_{\eye  j}^{-1}\pairing{g_\varepsilon' (h_n)c_{n}h_{n,t}}{\phi_\eye}$.
The vector $\dot{\eta}_n$ can be replaced by the function
$\Phi_n\left(\eta_n,\gamma_n\right)$ and thus we obtain a relation
\[
\frac{d\gamma_{n,j}}{dt}=\Psi_{n,j}\left(\eta_n,\gamma_n\right),
\]
 in place of Eq.~\eqref{eq:gamma_dot}, where $\Psi_{n,j}\left(\eta_n,\gamma_n\right)$
 is Lipschitz.  We therefore have a system of Lipschitz-continuous equations
\[
\frac{d\eta_{n,j}}{dt}=\Phi_{n,j}\left(\eta_n,\gamma_n\right),\qquad\frac{d\gamma_{n,j}}{dt}=\Psi_{n,j}\left(\eta_n,\gamma_n\right),
\]
and thus local existence theory~\cite{DoeringGibbon} guarantees
a solution
for the $\eta_{n,i}$'s
and $\gamma_{n,i}$'s for all times $t$ in a finite interval $0<t<\sigma$.
 This solution is, moreover, unique and continuous.  To continue this approximate
 solution to the PDE problem in Eq.~\eqref{eq:model_proof}
 up to an an arbitrary time $T_0>0$, it is necessary to find \emph{a priori}
 bounds on the approximate local-in-time solution.

\subsection{\emph{A priori} bounds on the Galerkin approximation}
\label{sec:a_priori_bds}

\noindent We identify the free energy
\[
F = \int_{\Omega}{dx}\,
\left[\tfrac{3}{2}h_x^2+G_\varepsilon(h)\right]
+ \int_{\Omega}{dx}\,
g_\varepsilon(h) \left[\tfrac{1}{4} \left(c^2-1\right)^2+\tfrac{1}{2}c_x^2\right],
\]
where $G_{\varepsilon}''(s)=1/\left[f_\varepsilon(s)g_\varepsilon(s)\right]$.
 Since the Galerkin approximation satisfies the weak form of the PDEs given
 in Eqs.~\eqref{eq:weak_h} and~\eqref{eq:weak_c}, it is possible to obtain
 the free-energy decay law
\begin{multline}
\frac{dF}{dt}(t)=
-\int_{\Omega-\Omega_-}{dx}\,  f_\varepsilon(h_n)
   \left[-h_{n,xxx}+\frac{h_{n,x}}{g_\varepsilon(h_n)
       h_\varepsilon(h_n)} + \frac{1}{g_\varepsilon
        (h_n)} \left(g_\varepsilon  (h_n)
       c_{n,x}^2\right)_x \right]^2dx \\- \int_\Omega{dx}\,
   g_\varepsilon  (h_n) \mu_{\varepsilon,n,x}^2%\\
+ \int_{\Omega_-}{dx}\,[\dots],\qquad
   0\leq t<\sigma,
\label{eq:free_energy_decay}
\end{multline}
where $\Omega_-(t)=\{x\in\Omega|h_n(x,t)<0\}$.
Now given the time-continuity of $h_n(x,t)$ in
$(0,\sigma)$, and the initial condition $h_n^0(x)>0$ (since
$n>n_0$), there is a time $\sigma_1>0$ such that
$h_n(x,t)>0$ for all $x\in\Omega$ and all
$t\in\left(0,\sigma_1\right)$.  Therefore,
$\Omega_-(t)=\emptyset$ for $t\in\left(0,\sigma_1\right)$, the last
integral in~\eqref{eq:free_energy_decay} vanishes, and hence
\[
F\left[c_n(x,t),h_n(x,t)\right] \leq
F\left[c_n(x,0),h_n(x,0)\right] \leq
\sup_{\varepsilon,n}F\left[c_n(x,0),h_n(x,0)\right]<\infty,
\]
for $0< t<\sigma_1$.  Consequently, we obtain the bound $\|h_{n,x}\|_2\leq
k_1$, where $0<t<\sigma_1$, and where $k_1$ depends only on the initial
conditions.  We have Poincar\'e's inequality for $h_{n,x}$,
\[
\|h_n\|_2^2\leq\left[\int_\Omega {dx}\, h_n(x)\right]^2+\left(\frac{L}{2\pi}\right)^2\|h_{n,x}\|_2^2.
\]
Now $\int_\Omega{dx}h_n(x,t)=L\eta_{n,0}(t)$.  Inspection
of Eq.~\eqref{eq:weak_h} shows that $\eta_{n,0}(t)=\eta_{n,0}\left(0\right)=\eta_{0}^0$.
 Thus,
\[
\|h_n\|_2^2\leq
L^2|\eta_{0}^0|^2+\left(\frac{L}{2\pi}\right)^2k_1\equiv k_2<\infty.
\] 
Using result~\eqref{eq:result1} from Appendix~\ref{apx:A},
we obtain the bound
\[
\|h_n\|_\infty\leq\frac{1}{\sqrt{L}}\|h_n\|_2+\sqrt{L}\|h_{n,x}\|_2
\equiv k_3\,.
\]
Additionally, the following properties hold:
\begin{itemize}
\item The function $h_{n}$ is H\"older continuous in space, with exponent
$\tfrac{1}{2}$\,;
\item $\int_\Omega {dx}\, G_{\varepsilon} (h_n)\leq k_4$\,.
\end{itemize}
These results hold in $0<t<\sigma_1$, and the constants $k_1$,
$k_2$, $k_3$, and $k_4$ are independent of $\varepsilon$, $n$, $\sigma$,
and $\sigma_1$, and in fact depend only on the functions $h_0(x)$
and $c_0(x)$.

To continue the estimates to the whole interval
$\left(0,\sigma\right)$, we need to prove that
$h_n\left(\cdot,\sigma_1\right)>0$ almost everywhere (a.e.).  If this
is true, there is a new interval $[\sigma_1,\sigma_2)$,
$\sigma_1<\sigma_2\leq\sigma$, on which $h_n(\cdot,t)>0$
a.e., and we can then provide \emph{a priori} bounds on
$h_n(\cdot,t)$ and $c_n(\cdot,t)$ on the
interval $\left[\sigma_1,\sigma_2\right)$.  It is then possible to
show that $h_n\left(\cdot.,\sigma_2\right)>0$ a.e.\ and thus, by
iteration, we extend the proof to the whole interval
$\left(0,\sigma\right)$, and find that $h_n\left(.,t\right)>0$ a.e.\
on $\left(0,\sigma\right)$.

We have the bound
\begin{equation}
\int_\Omega {dx}\, G_{\varepsilon}(h_n(\cdot,t))\leq k_4,
\label{eq:sigma_1_bound}
\end{equation}
where $k_4$ depends
only on the initial conditions, and where $0<t<\sigma_1$.  We now
specify $G_\varepsilon(s)$ in more detail.  This function
satisfies the condition
\[
G_\varepsilon''(s)=\frac{1}{f_\varepsilon(s)\,g_\varepsilon(s)}.
\]
We take $g_\varepsilon(s)$ to be as defined previously, and we choose a simple
regularisation for $f(s)$:
\[
f_\varepsilon(s)=g_\varepsilon(s)^3,
\]
which is Lipschitz continuous.
By defining
\[
\tilde{G}_\varepsilon(s)=-\int_s^\infty\frac{dr}{f_\varepsilon(r)g_\varepsilon(r)},\qquad
G_\varepsilon(s)=-\int_s^{\infty}{dr}\, \tilde{G}_\varepsilon(r),
\]
we obtain a function $G_\varepsilon(s)$ that is positive for all
$s\in\left(-\infty,\infty\right)$, and
\[
G_{\varepsilon}(s)=\tfrac{1}{6}\frac{1}{\left(s+\varepsilon\right)^2}\,,
\qquad s\geq0.
\]
Using the boundedness of $G_{\varepsilon}(s)$, and the time-continuity
of $h_n(\cdot,t)$, we apply the Dominated Convergence Theorem,
\[
\lim_{t\rightarrow\sigma_1}\int_\Omega {dx}\, G_{\varepsilon}(h_n(\cdot,t))=\int_\Omega{dx}\,
\lim_{t\rightarrow\sigma_1}G_{\varepsilon}(h_n(\cdot,t))=\int_{\Omega}{dx}\,
G_{\varepsilon}\left(h_n\left(\cdot,\sigma_1\right)\right)\leq
k_4.
\] 
Similarly, since the constant $k_1$ in the inequality $\|h_{n,x}\|_2\leq
k_1$, $0\leq t<\sigma_1$ depends
only on the initial data, we extend this last inequality to $t=\sigma_1$,
and thus $h_n(x,\sigma_1)$ is H\"older continuous in space.  

In the worst-case scenario, the time $\sigma_1$ is the first time at which
$h_n(x,t)$
touches down to zero, and and we therefore assume for contradiction that $h_n\left(x_0,\sigma_1\right)=0$,
and that $h_n(x,\sigma_1)\geq0$ elsewhere.
 Then, by H\"older continuity, for any $x\in\Omega$ we have the bound
 $0\leq h_n(x,\sigma_1)\leq k_1\left|x-x_0\right|^{\frac{1}{2}}$, and
 thus
\[
\int_\Omega {dx}\, G_{\varepsilon}\left(h_n\left(\cdot,\sigma_1\right)\right)
\geq\frac{k_1}{6}\int_0^L\frac{dx}{|x-x_0|+\varepsilon (2\sqrt{L}+\varepsilon)}.
\]
Hence,
\begin{multline*}
\frac{6}{k_1}\int_\Omega {dx}\, G_{\varepsilon}\left(h_n\left(\cdot,\sigma_1\right)\right)
\\
\geq-2\log\left[\varepsilon\left(2\sqrt{L}+\varepsilon\right)\right]+\log\bigg\{\left[L-x_0+\left(2\sqrt{L}+\varepsilon\right)\varepsilon\right]\left[x_0+\left(2\sqrt{L}+\varepsilon\right)\varepsilon\right]\bigg\}.
\end{multline*}
Thus, the integral $\int_\Omega G_\varepsilon(h_n(x,\sigma_1))dx$ can
be made arbitrarily large, which contradicts the
$\varepsilon$-independent bound for this quantity, obtained in
Eq.~\eqref{eq:sigma_1_bound}.  We therefore have the strong condition
that the set on which $h_n\left(\cdot,\sigma_1\right)\leq0$ is empty.
Iterating the argument, we have the important result
\begin{quote}
The set on which $h_n(\cdot,t)\leq0$ is empty, for $0<t<\sigma$.
\end{quote}
Using the same argument, we have an estimate on the minimum value of
$h_n(x,t)$,
\[
h_{\mathrm{min}}=\min_{x\in\Omega,t\in\left(0,\sigma\right]}h_n(x,t),
\]
namely,
\[
h_{\mathrm{min}}+\varepsilon\geq -k_1\sqrt{L}+\sqrt{k_1^2L+\frac{k_1^2L}{e^{k_4k_1^2}-1}},
\]
for all small positive $\varepsilon$.  Thus,
\begin{equation}
h_{\mathrm{min}}\geq M :=-k_1\sqrt{L}+\sqrt{k_1^2L+\frac{k_1^2L}{e^{k_4k_1^2}-1}},
\label{eq:min_h}
\end{equation}
a lower bound that depends only on the initial data $c_0(x)$ and $h_0(x)$.
 Note
that this result depends on the intimate relationship between the H\"older
continuity of a function and its boundedness in the $W^2$ norm, a relationship
that is true only in one dimension.  Thus, generalisation of this lower bound,
and by extension, long-time existence and uniqueness of solutions, does not
necessarily hold in higher dimensions.

Now, using Eq.~\eqref{eq:min_h} and the boundedness result 
\[
\int_\Omega{dx}\,
g_\varepsilon(h_n)
\left[\tfrac{1}{4}\left(c_n^2-1\right)^2+\tfrac{1}{2}c_{n,x}^2\right]
\leq k_5\,,
\]
where $k_5$ depends only on the initial data, we obtain an \emph{a priori}
bound on $\|c_{n,x}\|_2^2$,
\[
\int_\Omega{dx}\,  c_{n,x}^2 \leq \frac{2k_5}{M}\,.
\]
It is also possible to derive a bound on $\|c_n\|_2^2$.  We have the relation
\[
\int_\Omega{dx}\,\tfrac{1}{4}\left(c_n^2-1\right)^2\leq \frac{k_5}{M},
\]
which gives the inequality
$
\|c_n\|_4^4\leq 2\|c\|_2^2+\left({4k_5}/{M}\right).
$
Using the H\"older relation $\|c\|_2\leq |\Omega|^{\frac{1}{4}}\|c\|_4$,
we obtain a quadratic inequality in $\|c\|_2^2$,
\[
\|c\|_2^4\leq 2|\Omega|\|c\|_2^2+\frac{4|\Omega|k_5}{M},
\]
with solution
\[
\|c\|_2^2\leq |\Omega|+\frac{4|\Omega|k_5}{M},
\]
as required.  From the boundedness of $\|c_{n,x}\|_2$ and $\|c_{n}\|_2$ follows
the relation $\|c_n\|_\infty \leq k_6<\infty$, a result that depends only
on the initial conditions.
Let us recapitulate these results:
\begin{itemize}
\item $\|h_{n,x}\|_2$ is uniformly bounded;
\item $\|h_{n}\|_\infty$ is uniformly bounded;
\item The function $h_{n}$ is H\"older continuous in space, with exponent
$\tfrac{1}{2}$;
\item The function $h_{n}$ is nonzero everywhere and never decreases below
a certain value $M>0$, independent of $n$, $\varepsilon$, and $\sigma$.
\item $\|c_{n,x}\|_2$ is uniformly bounded;
\item $\|c_{n}\|_\infty$ is uniformly bounded;
\item The function $c_{n}$ is H\"older continuous in space, with exponent
$\tfrac{1}{2}$\,.
\end{itemize}
These results are independent of $n$, $\varepsilon$ and $\sigma$, and
hold for $0<t<\sigma$.
\subsection{Equicontinuity and convergence of the Galerkin approximation}
\label{sec:equicontinuity}

\noindent  Using Eq.~\eqref{eq:free_energy_decay}, we obtain the bound
\begin{multline*}
\int_0^t{dt'}\int_\Omega{dx}\, \bigg\{f_\varepsilon(h_n)\left[-h_{n,xxx}+\frac{h_{n,x}}{g_\varepsilon(h_n)f_\varepsilon(h_n)}+\frac{1}{g_\varepsilon(h_n)}\left(g_\varepsilon(h_n)c_{n,x}^2\right)_x\right]^2+g_\varepsilon(h_n)\mu_{\varepsilon,n,x}^2\bigg\}
\\
\leq F\left(0\right),\qquad 0<t<\sigma,
\end{multline*}
a bound that is independent of $n$, $\sigma$, and $\varepsilon$.  Since the
quantity $\|f_\varepsilon(h_n)\|_\infty=\left(\|h_n\|_\infty+\varepsilon\right)^3$
is bounded above by a constant $A_1$ independent of $n$, $\varepsilon$,
and $\sigma$, we have
\begin{multline*}
\int_0^tdt'\int_\Omega {dx}\,  J_{\varepsilon,n}^2\\
\leq A_1\int_0^t{dt'}\int_\Omega{dx}\, 
\bigg\{f_\varepsilon(h_n) \left[-h_{n,xxx} +
  \frac{h_{n,x}}{g_\varepsilon
     (h_n)f_\varepsilon(h_n)} +
  \frac{1}{g_\varepsilon(h_n)}
  \left(g_\varepsilon(h_n) c_{n,x}^2\right)_x\right]^2\bigg\}
\\
\leq A_1 F\left(0\right)\equiv A_2,
\end{multline*}
and thus
\[
\int_0^tdt'\,\|J_{\varepsilon,n}\|_2^2\leq A_2,\qquad 0< t<\sigma,
\]
where $A_2$ is independent of $n$, $\varepsilon$, and $\sigma$.  Similarly,
\[
\int_0^tdt'\,\|g_\varepsilon(h_n)\mu_{\varepsilon,n,x}\|_2^2\leq
A_3,\qquad 0<{t}<\sigma,
\]
and
\[
\int_0^tdt'\,\|c_n J_{\varepsilon,n}\|_2^2\leq A_2\|c_n\|_\infty^2\leq
A_2k_6^2,\qquad 0<t<\sigma,
\]
where $A_2$ and $A_3$ are independent of $n$, $\varepsilon$, and $\sigma$.
 By rewriting the evolution equations as
\[
h_{n,t}=-J_{\varepsilon,n,x},\qquad \left(g_\varepsilon(h)c\right)_t=-K_{\varepsilon,n,x},
\]
where $K_{\varepsilon,n}=c_n J_{\varepsilon,n}-g_\varepsilon(h_n)\mu_{\varepsilon,n,x}$,
we see that there are uniform bounds for $\int_0^tdt'\,\|J_{\varepsilon,n}\|_2^2$
and $\int_0^tdt'\,\|K_{\varepsilon,n}\|_2^2$, which depend only on the initial
data $c_0(x)$ and $h_0(x)$.

Bernis and Friedman~\cite{Friedman1990} proved the following claim:
\begin{quote}
\emph{Let $\varphi_i\left(x,t\right)$ be a sequence of functions, each of
which weakly satisfies the equation
\[
\varphi_{i,t}=-J_{i,x},\qquad J_i=J\left(\varphi_i\right).
\]
If $\varphi_i\left(x,\cdot\right)$ is H\"older continuous (exponent $\tfrac{1}{2}$),
and if the fluxes $J_i$ satisfy
\[
\int_0^t dt'\,\|J_i\|_2^2\leq B_1,\qquad 0<t<\sigma,
\]
where $B_1$ is a number independent of the index $i$ and the time $\sigma$,
then there is a constant $B_2$, independent of $i$ and $\sigma$, such that
\[
\left|\varphi_i\left(\cdot,t_2\right)-\varphi_i\left(\cdot,t_1\right)\right|\leq
B_2\left|t_2-t_1\right|^{\frac{1}{8}},
\]
for all $t_1$ and $t_2$ in $\left(0,\sigma\right)$.  
}
\end{quote}
We now observe that the fluxes $J_{\varepsilon,n}$, and $K_{\varepsilon,n}$
satisfy the conditions of this theorem, and thus
\begin{quote}
The functions $h_{n}(\cdot,t)$ and $c_n(\cdot,t)$
are H\"older continuous (exponent $\tfrac{1}{8}$), for $0<t<\sigma$.
\end{quote}
We therefore have a uniformly bounded and equicontinuous family of
functions $\{\left(h_n,c_n\right)\}_{n=n_0+1}^\infty$.  We also have a
recipe for constructing a uniformly bounded and equicontinuous
approximate solution
$\left(h_n(x,t),c_n(x,t)\right)$, in a small
interval $\left(0,\sigma\right)$.  The recipe can be iterated
step-by-step, and we obtain a uniformly bounded and equicontinuous
family of approximate solutions
$\{\left(h_n,c_n\right)\}_{n=n_0+1}^\infty$, on
$\left(0,T_0\right)\times\Omega$, for an arbitrary time $T_0$.  Then,
using the Arzel\`a--Ascoli theorem, we obtain the convergence result:
\begin{quote}
There is a subsequence of the family $\{\left(h_n,c_n\right)\}_{n=n_0+1}^\infty$
that converges uniformly to a limit $\left(h,c\right)$, in $\left[0,T_0\right]\times\Omega$.
\end{quote}
We prove several facts about the pair $\left(h,c\right)$.
\begin{quote}
\emph{Let $\left(h,c\right)$ be the limit of the family of functions
$\{\left(h_n,c_n\right)\}_{n=n_0+1}^\infty$
constructed in Secs.~\ref{sec:regularization}--\ref{sec:equicontinuity}.
Then the following properties hold for this limit:
\begin{enumerate}
\item The functions $h(x,t)$ and $c(x,t)$ are uniformly
H\"older continuous in space (exponent $\tfrac{1}{2}$), and uniformly H\"older
continuous in time (exponent $\tfrac{1}{8}$);
\item The initial condition $\left(h,c\right)\left(x,0\right)=\left(h_0,c_0\right)(x)$
holds;
\item $\left(h,c\right)$ satisfy the boundary conditions of the original
problem (periodic boundary conditions);
\item The derivatives $\left(h,c\right)_t$, $\left(h,c\right)_x$, $\left(h,c\right)_{xx}$,
$\left(h,c\right)_{xxx}$, and $\left(h,c\right)_{xxxx}$ are continuous in
space and time;
\item The function pair $\left(h,c\right)$ satisfy the weak form
of the PDEs,
\begin{multline*}
\int\int_{Q_{T_0}}dtdx\, h\varphi_t+\int\int_{Q_{T_0}}dtdx\,J_\varepsilon\varphi_x=0,\\
J_\varepsilon=f_\varepsilon(h)h_{xxx}-\frac{1}{g_\varepsilon(h)}h_x-\frac{f_\varepsilon(h)}{g_\varepsilon(h)}\left(g_\varepsilon(h)c_{x}^2\right)_x,
\end{multline*}
\begin{multline*}
\int\int_{Q_{T_0}}dtdx\, g_\varepsilon(h)c\psi_t+\int\int_{Q_{T_0}}dtdx\,K_\varepsilon\psi_x=0,\\
K_\varepsilon= cJ_\varepsilon-g_\varepsilon(h)\left[c^3-c-\frac{1}{g_\varepsilon(h)}\left(g_\varepsilon(h)c_x\right)_x\right],
\end{multline*}
where $\varphi\left(x,t\right)$ and $\psi\left(x,t\right)$ are suitable test
functions.
\label{page:claims45}
\end{enumerate}
}
\end{quote}
The statements~1,~2, and~3 are obvious.  Now, any pair
$\left(h_n(x,t),c_n(x,t)\right)$ satisfies the
equation set
\begin{multline*}
\int\int_{Q_{T_0}}dtdx\, h_n\phi_t+\int\int_{Q_{T_0}}dtdx\,J_{\varepsilon,n}\phi_x=0,
\\
J_{\varepsilon,n}=f_{\varepsilon} (h_n)h_{n,xxx}-\frac{1}{g_\varepsilon(h_n)}h_x-\frac{f_\varepsilon(h_n)}{g_\varepsilon(h_n)}\left(g_\varepsilon(h_n)c_{n,x}^2\right)_x,
\end{multline*}
\begin{multline*}
%\hskip -0.2in
\int\int_{Q_{T_0}}dtdx\, g_\varepsilon(h_n)c_n\psi_t+\int\int_{Q_{T_0}}dtdx\,
K_{\varepsilon,n}\psi_x=0,\\
K_{\varepsilon,n}= c_nJ_{\varepsilon,n}-g_\varepsilon
(h_n)\left[c_n^3-c_n-\frac{1}{g_\varepsilon(h_n)}\left(g_\varepsilon
    (h_n) c_{n,x}\right)_x\right],
\end{multline*}
and from the boundedness of the fluxes $J_{\varepsilon,n}$ and $K_{\varepsilon,n}$
in $L^2\left(0,T_0;L^2(\Omega)\right)$, it follows that
\[
\left(J_{\varepsilon,n}\,,\,K_{\varepsilon,n}\right)\rightharpoonup\left(J_\varepsilon\,,\,K_\varepsilon\right),
\qquad\text{weakly in }L^2\left(0,T_0;L^2(\Omega)\right),
\]
for a subsequence.  Using the regularity
theory for uniformly parabolic equations and the uniform H\"older continuity
of the $\left(h_n,c_n\right)$'s, it follows that
\begin{quote}
The derivatives $\left(h_n,c_n\right)_t$, $\left(h_n,c_n\right)_x$,
$\left(h_n,c_n\right)_{xx}$, $\left(h_n,c_n\right)_{xxx}$, and $\left(h_n,c_n\right)_{xxxx}$
are uniformly convergent in any compact subset of $\left(0,T_0\right]\times\Omega$.
\end{quote}
Thus,
\[
J_\varepsilon = f_\varepsilon(h)h_{xxx} -
\frac{1}{g_\varepsilon(h)} h_x -
\frac{f_\varepsilon(h)}{g_\varepsilon(h)} \left(g_\varepsilon(h)c_{x}^2\right)_x,
\]
\[
K_\varepsilon=cJ_\varepsilon-g_\varepsilon(h)\left[c^3-c-\frac{1}{g_\varepsilon(h)}\left(g_\varepsilon(h)c_x\right)_x\right],
\]
on $\left(0,T_0\right]\times\Omega$, and therefore, claims~4 and~5 on
p.~\pageref{page:claims45} follow.

\subsection{Convergence of regularised problem, as $\varepsilon\rightarrow0$}
\label{sec:convergence}

\noindent The result in Sec.~\ref{sec:equicontinuity} produced a solution
$\left(h_\varepsilon,c_\varepsilon\right)$ to the regularised problem.  Due
to the result
\begin{equation}
h_{\varepsilon}\left(x,t\right)\geq h_{\mathrm{min}}\geq M=-k_1\sqrt{L}+\sqrt{k_1^2L+\frac{k_1^2L}{e^{k_4k_1^2}-1}}>0,\\
\label{eq:no_rupture}
\end{equation}
independent of $\varepsilon$, the argument of
Sec.~\ref{sec:equicontinuity} can be recycled to produce a solution
$\left(h,c\right)$ to the unregularised problem.  This solution is
constructed as the limit of a convergent subsequence, formally written as
$\left(h,c\right)=\lim_{\varepsilon\rightarrow0}\left(h_\varepsilon,c_\varepsilon\right)$,
and the results of the theorem in Sec.~\ref{sec:equicontinuity}
apply again to $\left(h,c\right)$.  The result~\eqref{eq:no_rupture}
applies to $h$ constructed as
$h=\lim_{\varepsilon\rightarrow0}h_\varepsilon$, and thus all the
derivatives $\left(h,c\right)_t$, $\left(h,c\right)_x$,
$\left(h,c\right)_{xx}$, $\left(h,c\right)_{xxx}$, and
$\left(h,c\right)_{xxxx}$ are continuous on the whole space
$\left(0,T_0\right]\times\Omega$ and therefore, the weak solution
$\left(h,c\right)$ is in fact a strong one.
 
\subsection{Regularity properties of the solution $\left(h,c\right)$}
\label{sec:regularity}

\noindent 
Using a bootstrap argument, we show that the solution $\left(h,c\right)$ belongs to the regularity
classes $L^{\infty}\left(0,T_0;\Sobo^{2,2}(\Omega)\right)$ and $L^2\left(0,T_0;\Sobo^{4,2}(\Omega)\right)$.
 From Sec.~\ref{sec:a_priori_bds} it follows immediately that
\[
\|h_x\|_2\,,\ \|c_x\|_2 < \infty,
\]
with time-independent bounds.  Thus, using Poincar\'e's inequality, it follows
that $h,c\in \Sobo^{1,2}(\Omega)$ and, moreover, 
\[
\sup_{\left[0,T_0\right]}\|h_x\|_2\,,\ \sup_{\left[0,T_0\right]}\|c_x\|_2 <
\infty.
\]
From Sec.~\ref{sec:equicontinuity},
it follows that $J$ and $\mu_x$ belong to the regularity class $L^2\left(0,T_0;L^2(\Omega)\right)$,
and hence $J$, $\mu_x\in L^2\left(0,T_0;L^1(\Omega)\right)$. 
The functions $J$, $\mu$, and $\mu_x$ take the form
\[
J=h^3 h_{xxx}-h_xh^{-1}-h^2\left(h_x c_x^2+2h c_x c_{xx}\right),
\]
\[
\mu = c^3-c - h^{-1}\left(h_x c_x + hc_{xx}\right),
\]
and
\[
\mu_x= \left(3c^2-1\right)c_x + h^{-2}h_x^2 c_x - h^{-1}h_{xx}c_x - h^{-1}h_xc_{xx}
-c_{xxx},
\]
respectively.
We make use of the following observations:
%The following results will help us in our demonstration,
%
%
\begin{itemize}
\item The function $h(x,t)$ is bounded from above and below,
\[
0< h_{\mathrm{min}}\leq h(x,t)\leq h_{\mathrm{max}}<\infty;
\]
we also have the boundedness of $c(x,t)$, $\|c\|_\infty
<\infty$. 
\item Since $\mu_x\in L^2\left(0,T_0;L^2(\Omega)\right)$, it follows
that $\mu\in L^2\left(0,T_0; L^2(\Omega)\right)$, by Poincar\'e's
inequality.
\item From this it follows that $c_x h_x + hc_{xx}$ is in the class $L^2\left(0,T_0;L^2(\Omega)\right)$.
\item  Given the inequality 
\[
\|h\mu c_x\|_2^2\leq h_{\mathrm{max}}\|\mu\|_\infty^2\|c_x\|_2^2\leq
h_{\mathrm{max}}\|c_x\|_2^2\left[\frac{1}{\sqrt{L}}\|\mu\|_2+\sqrt{L}\|\mu_x\|_2\right]^2,
\]
we have the result $h_x c_x^2 + h c_x c_{xx}\in L^2\left(0,T_0;L^2(\Omega)\right)$.
\item Similarly, since $\int_0^{T_0}\|h\mu h_x\|_2^2dt<\infty$, we have the
bound
$h_x^2c_x +hh_x c_{xx}\in  L^2\left(0,T_0;L^2(\Omega)\right)$.
\end{itemize}
Using these facts, and relegating the details to Appendix~\ref{apx:B},
it is possible to show that $c_{xx}$ is in
$L^2\left(0,T_0;L^1(\Omega)\right)$, from which follows the result
$h_{xxx},c_{xxx}\in L^1\left(0,T_0;L^1(\Omega)\right)$.  These results
give rise to further bounds, namely $c_{xx}\in
L^2\left(0,T_0;L^2(\Omega)\right)$, and
$\int_0^{T_0}dt\,\|h_x^2c_x\|_2<\infty$, whence $h_{xx}\in
L^2\left(0,T_0;L^2(\Omega)\right)$.  Using this collection of bounds,
we obtain
\[
h,c\in L^1\left(0,T_0;\Sobo^{3,1}(\Omega)\right),
\]
and hence finally, 
\[
h,c\in L^2\left(0,T_0;\Sobo^{3,2}(\Omega)\right).
\]
Thus, the solution $\left(h,c\right)$ belongs to the following regularity
class:
\begin{equation}
\left(h,c\right)\in L^{\infty}\left(0,T_0;\Sobo^{1,2}(\Omega)\right)
\cap L^2\left(0,T_0;\Sobo^{3,2}(\Omega)\right)
\cap C^{\frac{1}{2},\frac{1}{8}}\left(\Omega\times\left[0,T_0\right]\right).
\label{eq:regularity0}
\end{equation}

Extra regularity is obtained by writing the equation pair as
\begin{align*}
\frac{\partial h}{\partial t}+h^3h_{xxxx}&=-3h^2h_xh_{xxx}+\varphi_1+\varphi_2\equiv
\varphi\left(x,t\right),\\
\frac{\partial c}{\partial t}+c_{xxxx}&=-\frac{2}{h}h_xc_{xxx}+\psi_1+\psi_2\equiv
\psi\left(x,t\right),
\end{align*}
where 
\begin{align*}
\int_0^{\tau}{dt}\,\|\varphi_1\|_2&\leq \sup_{\left[0,\tau\right]}\|c_{xx}\|_2\int_0^{\tau}{dt}\,|\nu_1|,\qquad
\nu_1\in L^2\left(\left[0,\tau\right]\right),\\
\int_0^{\tau}{dt}\,\|\psi_1\|_2&\leq \sup_{\left[0,\tau\right]}\|c_{xx}\|_2\int_0^{\tau}{dt}\,|\nu_2|,\qquad
\nu_2\in L^2\left(\left[0,\tau\right]\right),
\end{align*}
for any $\tau\in \left(0,T_0\right]$, and where $\varphi_2$ and $\psi_2$
belong to
the class $L^2\left(0,T_0;L^2(\Omega)\right)$.  Here we have used
the form of the concentration equation given by Eq.~\eqref{eq:conc_no_flux}.
 By multiplying
the height and concentration equations by $h_{xxxx}$ and $c_{xxxx}$ respectively,
and by integrating over space and time, it is readily shown that
\[
\left(h,c\right)\in L^\infty\left(0,T_0;\Sobo^{2,2}(\Omega)\right),
\]
and hence
\[
\left(\varphi,\psi\right)\in L^2\left(0,T_0;L^{2}(\Omega)\right),
\]
from which follows the regularity result
\begin{equation}
\left(h,c\right)\in L^{\infty}\left(0,T_0;\Sobo^{2,2}(\Omega)\right)
\cap L^2\left(0,T_0;\Sobo^{4,2}(\Omega)\right)
\cap C^{\frac{3}{2},\frac{1}{8}}\left(\Omega\times\left[0,T_0\right]\right).
\label{eq:regularity}
\end{equation}

\subsection{Uniqueness of solutions}
\label{sec:uniqueness}
Let us consider two solution pairs $\left(h,c\right)$ and $\left(h',c'\right)$
and form the difference $\left(\delta h,\delta c\right)=\left(h-h',c-c'\right)$.
 Given the initial conditions $\left(\delta c(x,0),\delta h(x,0)\right)=\left(0,0\right)$,
 we show that $\left(\delta h,\delta c \right)=\left(0,0\right)$ for all
 time, that is, that the solution we have constructed is unique.
We observe that that the equation for the difference $\delta c$ can be written
in the form
\begin{equation}
\frac{\partial}{\partial t}\delta{c}+\frac{\partial^4}{\partial{x}^4}\delta{c}=\delta
\varphi\left(x,t\right),
\label{eq:delta_c_unique}
\end{equation}
where 
$\delta \varphi\left(x,t\right) \in L^2\left(0,T_0;L^2(\Omega)\right)$,
and where $\delta \varphi\left(\delta c=0\right)=0$.
Using semigroup theory~\cite{TaylorPDEbook},
we find that Eq.~\eqref{eq:delta_c_unique} has a unique solution.
Since $\delta{c}=0$ satisfies Eq.~\eqref{eq:delta_c_unique}, and since $\delta{c}\left(x,0\right)=0$,
it follows that $\delta{c}=0$ for all times $t\in\left[0,T_0\right]$.

It is now possible to formulate an equation for the difference $\delta{h}$
by subtracting the evolution equations of $h$ and $h'$ from one another,
mindful that $\delta c=0$.  We multiply the resulting equation by
$\delta{h}_{xx}$, integrate over space, and obtain after using
inequalities (see Appendix~\ref{apx:C})
\begin{multline*}
2\kappa\sup_{\tau\in\left[0,T\right]}\|\delta{h}_x\|_2^2\left(\tau\right)
\leq\sup_{\tau\in\left[0,T\right]}\|\delta{h}_x\|_2^2\left(\tau\right)\int_0^T{dt}\,\|h^{-1}+h'^2c_x^2\|_\infty^2\\
+\kappa_P^2\sup_{\tau\in\left[0,T\right]}\|\delta{h}_x\|_2^2\left(\tau\right)\int_0^T{dt}\left(\|h'_{xxx}+4c_xc_{xx}\|_2^2+\|h_x'+h_xc_x^2\|_2^2\right),
\end{multline*}
where $\kappa$ and $\kappa_P$ are numerical constants.
Using the results of Sec.~\ref{sec:regularity}, it is
readily shown that $h^{-1}+h'^2c_x^2\in L^2\left(0,T;L^{\infty}(\Omega)\right)$,
and that the functions $h'_{xxx}+4c_xc_{xx}$ and $h_x'+h_xc_x^2$ belong to
the class $L^2\left(0,T;L^2(\Omega)\right)$.  By choosing $T$
sufficiently small, it is possible to impose the inequality
\[
\frac{1}{2\kappa}\left[\int_0^T{dt}\,\|h^{-1}+h'^2c_x^2\|_\infty^2+\kappa_P^2\int_0^T{dt}\left(\|h'_{xxx}+4c_xc_{xx}\|_2^2+\|h_x'+h_xc_x^2\|_2^2\right)\right]<1,
\]
which in turn forces $\sup_{\tau\in\left[0,T\right]}\|\delta{h}_x\|_2^2=0$,
and hence the solution is unique.

\section{Parametric dependence of the height dip}
\label{sec:height_dip}

\noindent 
In this section we perform numerical simulations of the equations~\eqref{eq:model}
and focus on one feature of the equations: the drop in the free-surface height
at the boundary between binary fluid domains.  Our results in Sec.~\ref{sec:existence}
gave a rigorous upper bound for the magnitude of this height dip, as a function
of the problem parameters, and we compare this estimate with numerical solutions.

\begin{figure}[htb]
\subfigure[]{
 \includegraphics[width=.45\textwidth]{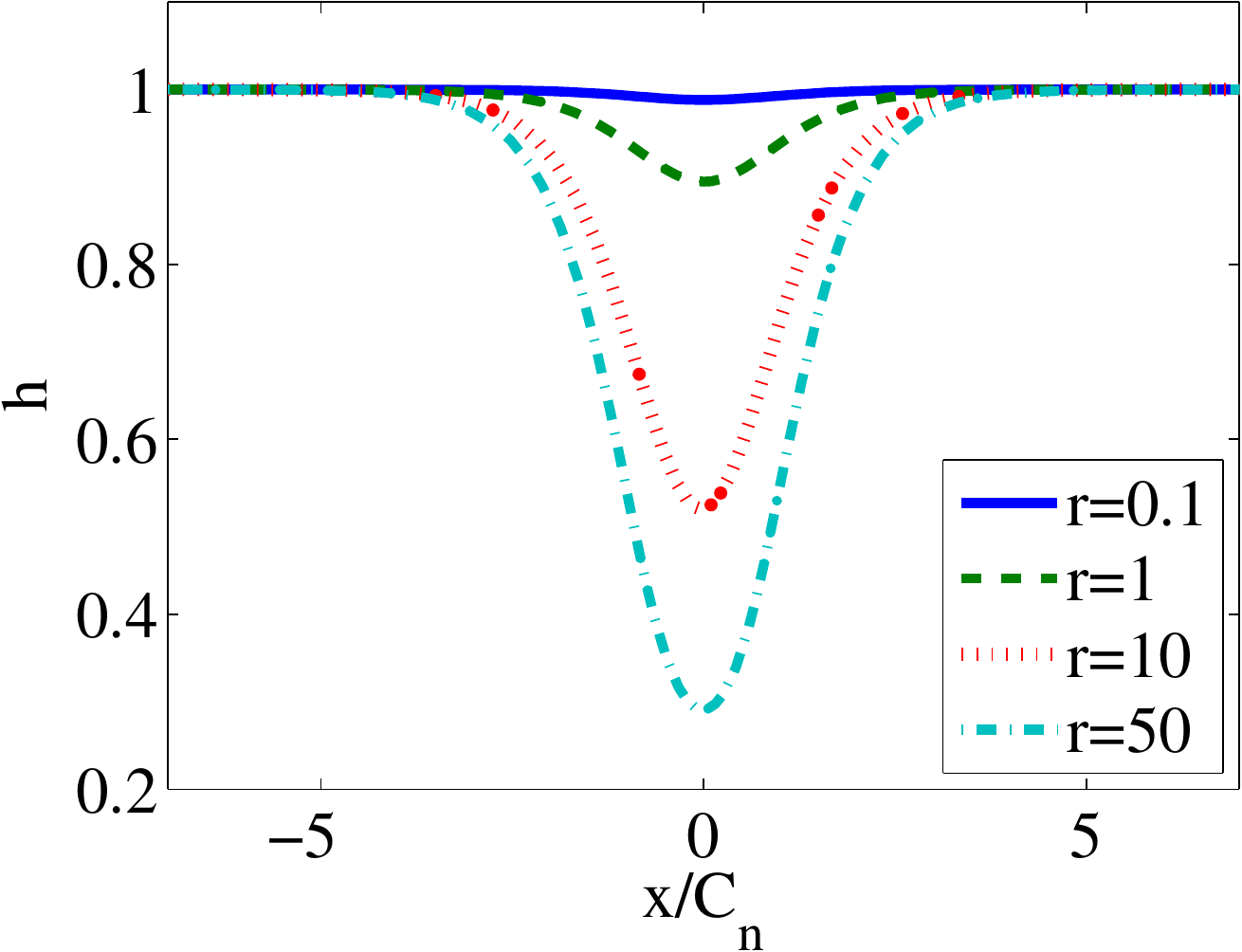}
}\hspace{1em}
\subfigure[]{
 \includegraphics[width=.45\textwidth]{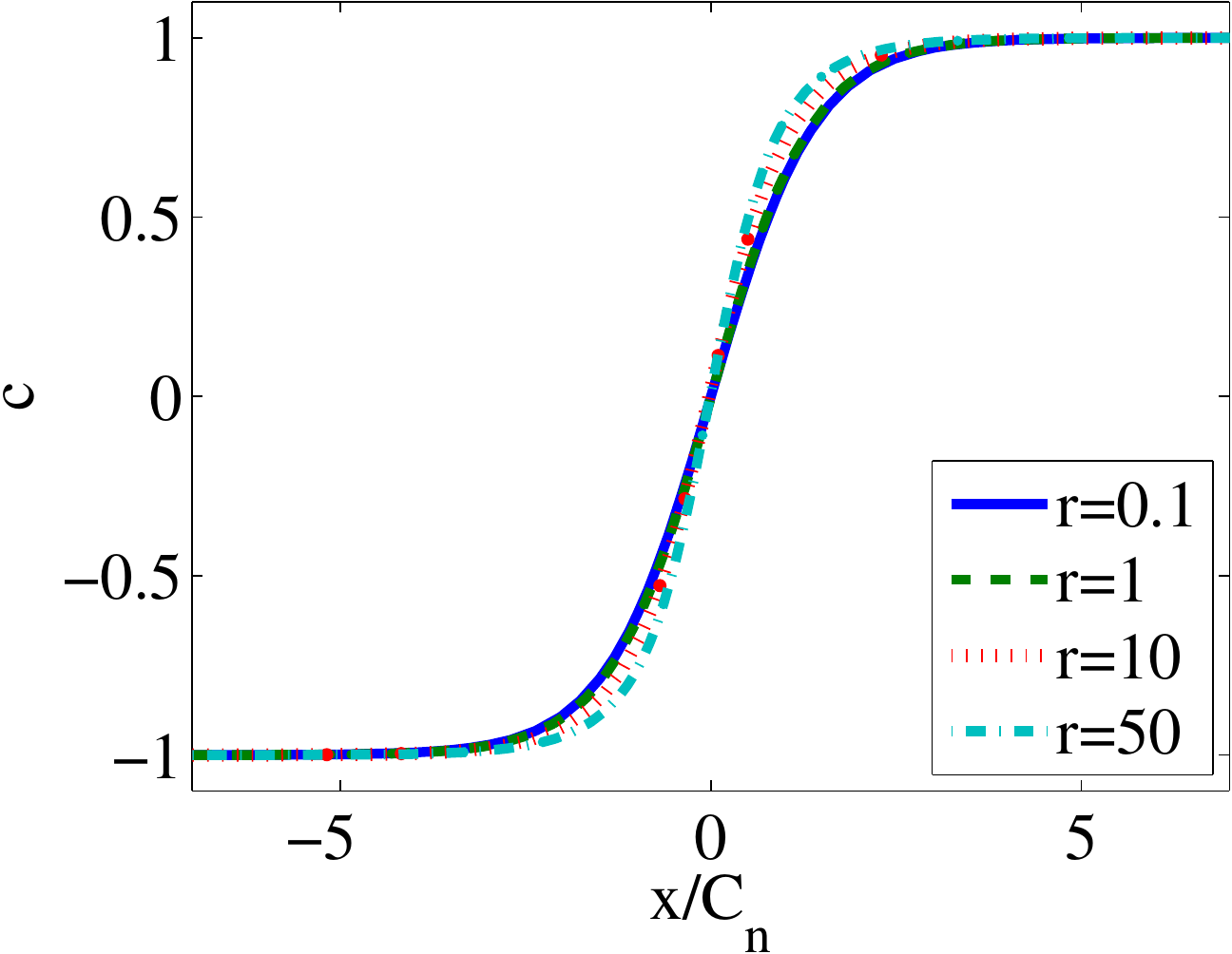}
}
\caption{%(Color online)  
Equilibrium solutions of Eq.~\eqref{eq:eqm}
for $C=C_{\mathrm{n}}^2|A|=1$ and $r=0.1,1,10,50$.  In (a) the valley deepens
with increasing $r$ although the film never ruptures, while in (b) the front
steepens with increasing $r$. (From \'O N\'araigh and
Thiffeault~\cite{ONaraigh2007}.)
}
\label{fig:hc}
\end{figure}

We perform numerical simulations of the full equations~\eqref{eq:model},
with initial data comprising a perturbation away from the unstable steady
state $\left(h,c\right)=\left(1,0\right)$.
The free surface and concentration evolve to an equilibrium state where the
salient feature is the formation of domains (intervals where $c\approx\pm1$)
that are separated by smooth transition regions, across which the free surface
dips below its average value.  Since we are interested in this characteristic
asymptotic feature of phase separation, we shift our focus instead to the
steady version from Eq.~\eqref{eq:model} obtained by setting $\partial_t=J=\mu_x=0$.
 For this reason also, we consider only the repulsive Van der Waals force,
 for
 which the interaction potential is given by $\phi=-|A|h^{-3}$.
We then solve the boundary-value problem
\begin{subequations}
\begin{equation}
\frac{1}{C}\frac{\partial^2 h}{\partial x^2}=|A| C_{\mathrm{n}}^2\left(1-\frac{1}{h^3}\right)
+r\left[\tfrac{1}{4}\left(c^2-1\right)^2+\tfrac{1}{2}\left(\frac{\partial
c}{\partial x}\right)^2\right],
\end{equation}
\begin{equation}
\frac{\partial^2 c}{\partial x^2}=c^3-c-\frac{1}{h}\frac{\partial
h}{\partial
x}\frac{\partial c}{\partial x},
\end{equation}%
\label{eq:eqm}%
\end{subequations}%
where now the domain is infinite and the boundary conditions are $h\left(\pm\infty\right)=1$,
and
$\mu\left(\pm\infty\right)=0$.  We have rescaled lengths by $C_{\mathrm{n}}$.
In Fig.~\ref{fig:hc} we present numerical solutions
exhibiting the dependence of the solutions on the parameters in Eq.~\eqref{eq:eqm}.
As before,
the free-surface height possesses peaks and valleys, where the valleys occur
in the concentration field's transition zone.  These profiles are qualitatively
similar to the results obtained in experiments on thin binary films~\cite{WangW2003,
WangH2000, ChungH2004}.  While the valley increases in depth
\begin{figure}[htb]
\includegraphics[width=.5\textwidth]{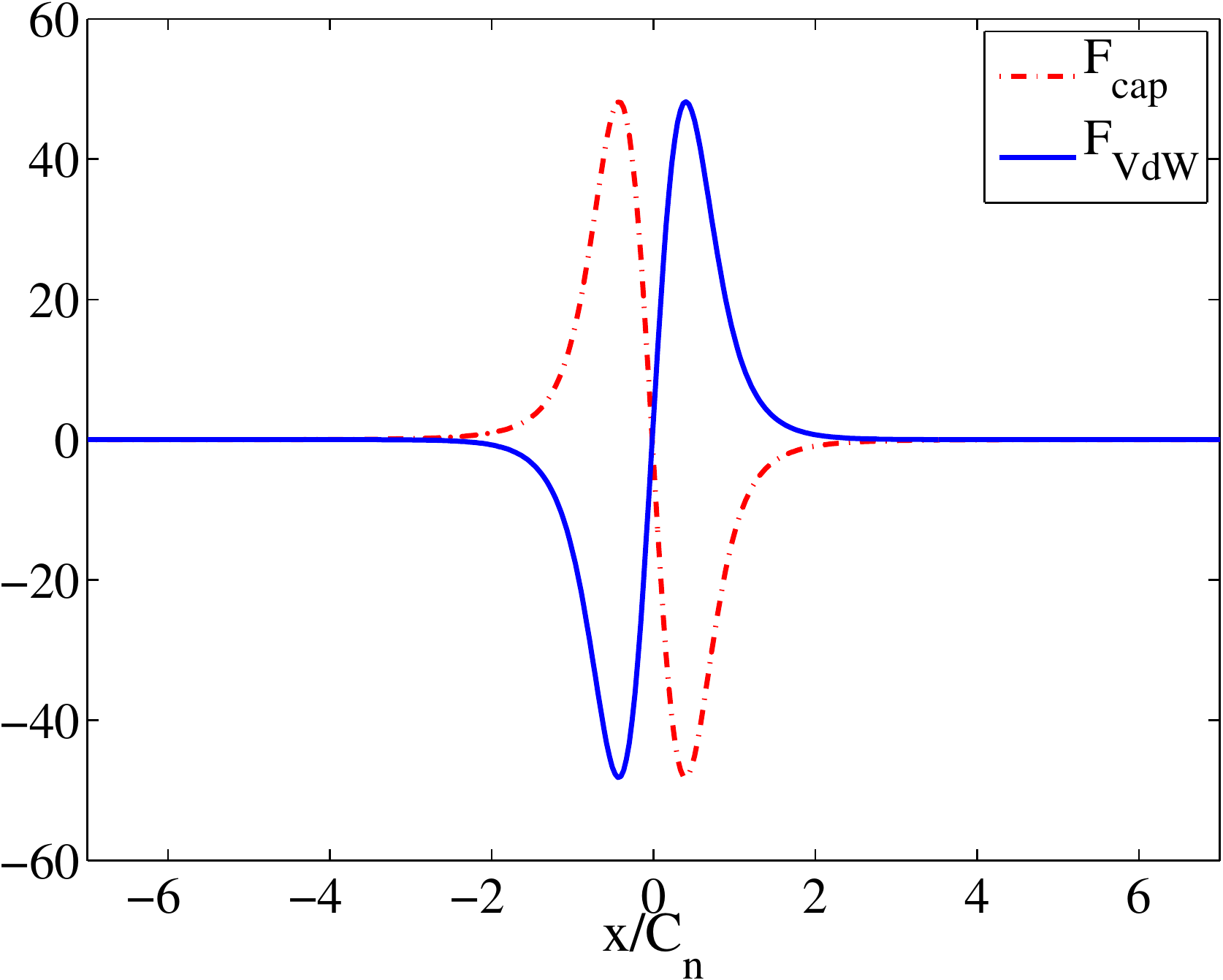}
\caption{The forces $F_{\mathrm{cap}}$ and $F_{\mathrm{VdW}}$
for $C=C_{\mathrm{n}}^2=|A|=1$ and $r=50$ are shown to have opposite sign.}
\label{fig:vdw_cap}
\end{figure}
for large $r$, rupture never takes place, in agreement with the theory in
Sec.~\ref{sec:existence}.  Similar results are obtained by varying $|A|$, with
smaller $|A|$ leading to larger dips.  This behaviour is underlined by the
results in Fig.~\ref{fig:vdw_cap}, where we plot the sign and magnitude of
the backreaction or capillary force $F_{\mathrm{cap}}=-rh^{-1}\left(hc_x^2\right)_x$
and the Van der Waals force $F_{\mathrm{VdW}}=|A|\left(h^{-3}\right)_x$.
 These forces have opposite sign: the Van der Waals force inhibits film rupture,
 while the backreaction promotes film thinning.    The thinning of binary
 films due to capillary or backreaction effects has been documented in experiments~\cite{WangH2000}.
 
\begin{figure}[htb]
        \subfigure[]{
        \includegraphics[width=.45\textwidth]{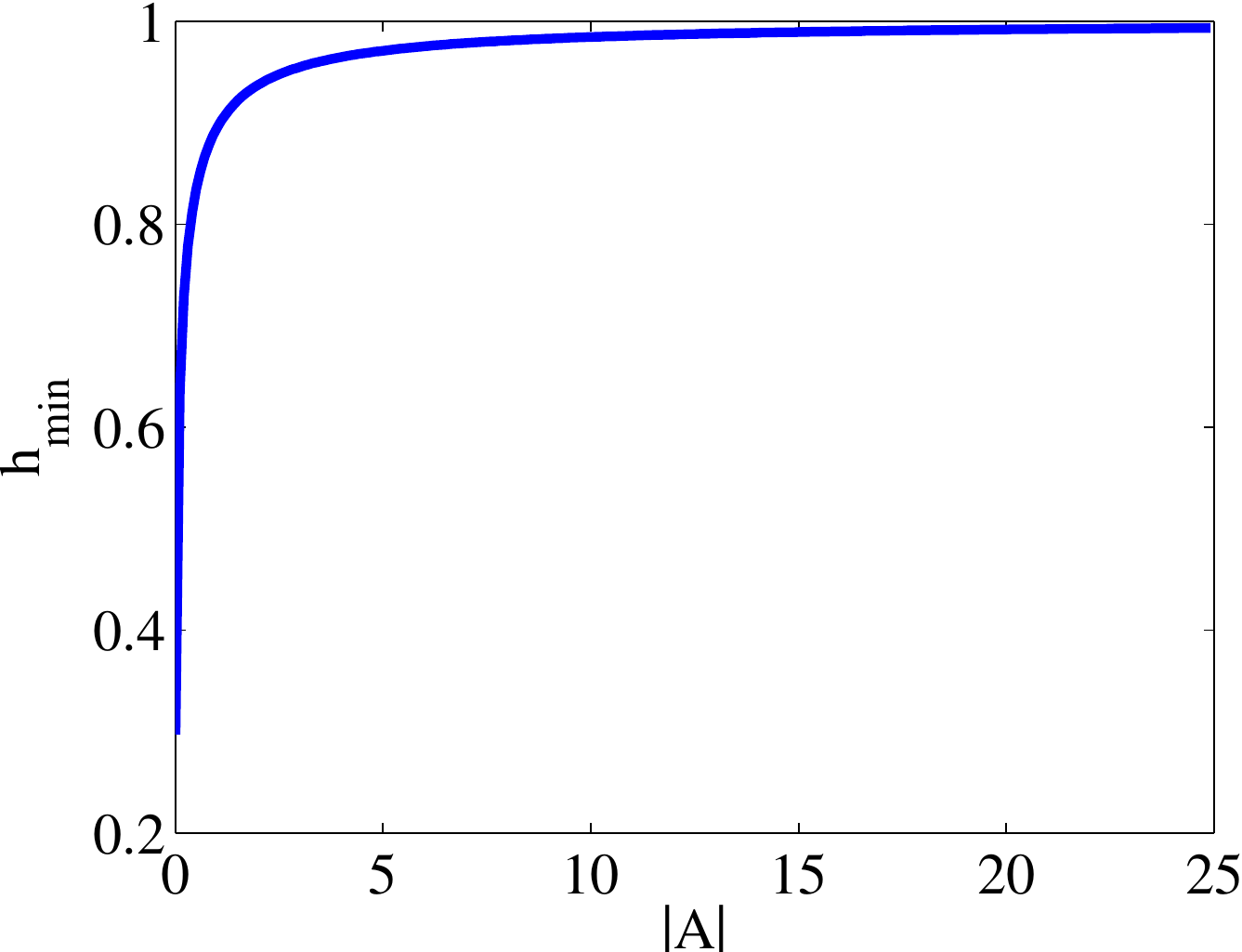}
}\hspace{1em}
        \subfigure[]{
        \includegraphics[width=.45\textwidth]{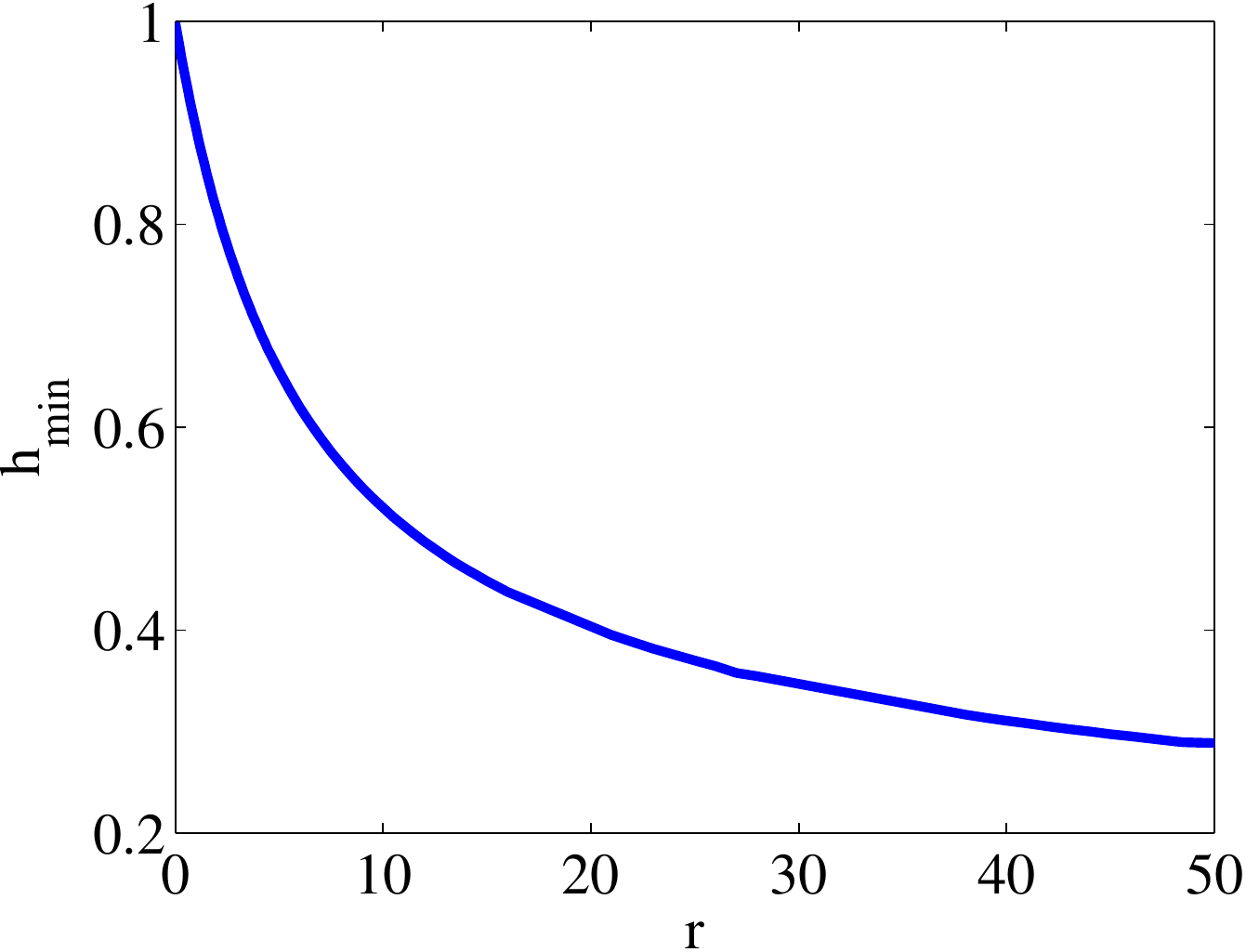}
}
\caption{Dependence of dip magnitude $h_{\mathrm{min}}$ on (a) the parameter
$|A|$; (b) the parameter $r$.  We have set $C=C_{\mathrm{n}}=1$.  In (a) the
dip decreases in magnitude with increasing $|A|$, demonstrating the tendency
of the Van der Waals force to flatten the free surface, while in (b) the
dip increases with increasing $r$, which shows how the backreaction thins
the film across transition zones.}
\label{fig:scaling_r_A}
\end{figure}

\begin{figure}[htb]
\includegraphics[width=.5\textwidth]{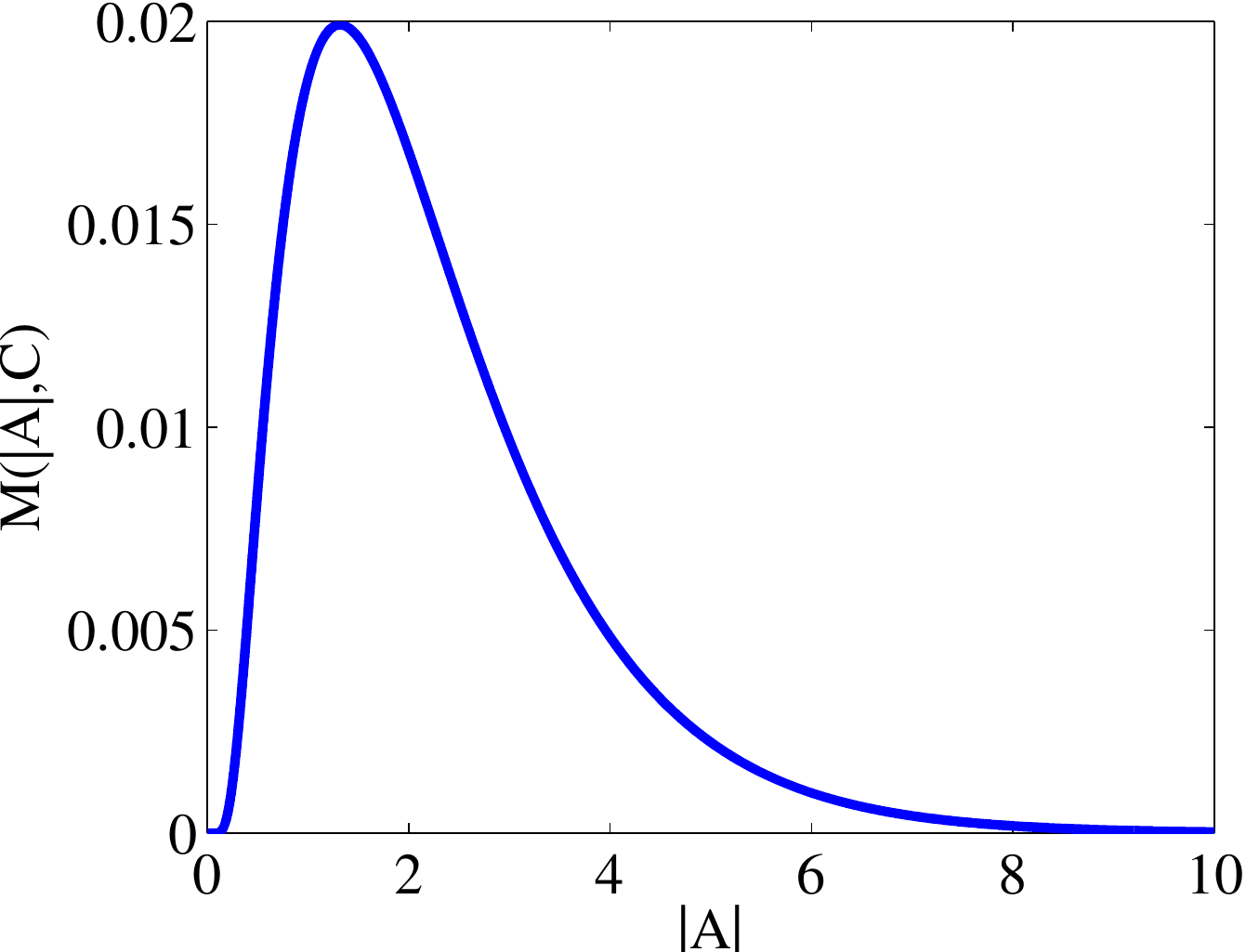}
\caption{A typical plot of $M(|A|,C)$ for $F_0=F_1=\tfrac{1}{2}$
and $C=1$.
 This theoretical lower bound has a different shape from those in Fig.~\ref{fig:scaling_r_A},
 which suggests that while $M(|A|,C)$ plays an important role in
 the analysis of the model equations, it does not capture the physics of
 film thinning.}
\label{fig:M_A}
\end{figure}

We are interested in the magnitude of the dip in the free-surface height,
as a function of the problem parameters, and we plot the dependence of $h_{\mathrm{min}}$
as a function of $|A|$ and $r$ in Fig.~\ref{fig:scaling_r_A}.
We compare these results with the estimate for $h_{\mathrm{min}}$ found in
Sec.~\ref{sec:existence}.  In terms of the physical parameters of the
system, this estimate is
\[
h_{\mathrm{min}}\geq M\equiv\sqrt{2CL(F_0+F_1|A|)}\left(\sqrt{\frac{e^{4C|A|^{-1}\left(F_0+F_1|A|\right)^2}}{e^{4C|A|^{-1}\left(F_0+F_1|A|\right)^2}-1}}-1\right),
\]
where $F_1=\tfrac{1}{2}\int_\Omega{dx}\,\left[h(x,0)\right]^{-2}\neq0$, and
$F_0=F\left(0\right)-F_1$. 
The function $M(|A|,C)$ has no explicit $r$-dependence: although
$F_0$ depends on $r$, it is possible to find initial data to remove this
dependence.
We show a representative plot of $M(|A|,C)$ in Fig.~\ref{fig:M_A}.
Although a comparison between Fig.~\ref{fig:scaling_r_A} and Fig.~\ref{fig:M_A}
is not exact, since the boundary conditions and domains are different in
both cases, we see that the shape of the bound in Fig.~\ref{fig:scaling_r_A}
is different from that in Fig.~\ref{fig:M_A}.  Since the bound in Fig.~\ref{fig:scaling_r_A}
is obtained from numerical simulations, and is intuitively correct, we conclude
that it has the correct shape and that the bound of Fig.~\ref{fig:M_A}, while
mathematically indispensable, is not sharp enough to be useful in determining
the parametric dependence of the dip in free-surface height.

\section{Conclusions}
\label{sec:conclusions}

\noindent Starting from the Navier--Stokes Cahn--Hilliard equations, we have
derived a pair of nonlinear parabolic PDEs that describe the coupled 
effects of phase separation and free-surface variations in a thin film of
binary liquid.  Since we are interested in the long-time outcome of the phase
separation, we focused on liquids that experience a repulsive Van der Waals
force, which tends to inhibit film rupture.  Using physical intuition, we
identified a decaying energy functional that facilitated analysis of the
equations.
We have shown that given sufficiently smooth initial data, solutions to the
model equations~\eqref{eq:model} exist in a strong sense, and
have certain regularity properties.  Our proof follows the method developed
in~\cite{Friedman1990}.  Central to the analysis is the decaying free energy,
and the derivation of a no-rupture condition, which prevents the film from
touching down to zero.  The no-rupture condition is valid in one spatial
dimension only, and thus existence and uniqueness results will not necessarily
be obtainable in higher dimensions.

We carried out one-dimensional numerical simulations of the full
equations and found that the free-surface height and concentration
tend to an equilibrium state.  The concentration forms domains; that
is, extended regions where $c\approx\pm1$.  The domains are separated
by narrow zones where the concentration smoothly transitions between
the limiting values $\pm1$.  At the transition zones, the free surface
dips below its mean value, a feature of binary thin-film behaviour
that is observed in experiments.  To study the magnitude of this dip
as a function of the problem parameters, we focused on solving the
equilibrium version of Eq.~\eqref{eq:model} as a boundary-value
problem.  This simplification is carried out without loss of
generality, since we have shown that the system tends asymptotically
to such a state.  We have shown that the magnitude of the dip
decreases by increasing the strength of the repulsive Van der Waals
force, while the dip depth actually increases by increasing the
strength of the backreaction.  Thus, in the absence of the Van der
Waals force, the film would rupture, preventing the occurrence of the
phase separation so characteristic of long-time Cahn--Hilliard
dynamics.  The film-thinning tendency of the backreaction has been
observed experimentally~\cite{WangH2000, WangW2003, ChungH2004}.
Simulations involving two lateral directions have elsewhere been carried
out by the authors~\cite{ONaraigh2007}, and the qualitative features are
similar to those obtained here.

\emph{Acknowledgements.} We thank G.~A.~Pavliotis for helpful discussions.
L.\'O.N. was supported by the Irish government and the UK Engineering
and Physical Sciences Research Council. J.-L.T. was supported in part
by the UK EPSRC Grant No. GR/S72931/01.

\appendix

\section{}
\label{apx:A}
\noindent In this appendix we give a list of nonstandard inequalities used in
the paper and provide proofs.

\begin{enumerate}
\item Let $\phi:\Omega\subset\mathbb{R}\rightarrow\mathbb{R}$ belong
to the class $\Sobo^{1,2}(\Omega)$.  Then the following string of inequalities
holds,
%the following inequality holds,
%
%
%
\begin{equation}
\sup_\Omega |\phi|\leq \frac{1}{L}\|\phi\|_1 + \|\phi_x\|_1
\leq\frac{1}{\sqrt{L}}\,\|\phi\|_2+\sqrt{L}\,\|\phi_x\|_2,
\label{eq:result1}
\end{equation}
\emph{Proof:} Using the Fundamental Theorem of Calculus, we have
\[
\phi\left(x\right)=\phi\left(a\right)
+ \int_a^x{ds}\,\frac{\partial\phi}{\partial s},
\]
for any distinct points $x$ and $a$ in $\Omega$.
Since the function $|\phi(x)|$ is continuous on $\Omega$, it has a maximum
value in $\Omega$, attained at the point $x_{\mathrm{max}}$.
Thus,
\[
|\phi\left(x_{\mathrm{max}}\right)|\equiv\sup_{\Omega}|\phi|\leq |\phi\left(a\right)|+\int_a^{x_{\mathrm{max}}}{ds}\left|\frac{\partial\phi}{\partial
s}\right|\leq  |\phi\left(a\right)|+\int_\Omega{ds}\left|\frac{\partial\phi}{\partial
s}\right|.
\]
Since this is true for any $a\in\Omega$, by integrating over $a$, we obtain
the inequality
\[
\sup_\Omega|\phi|\leq\frac{1}{L}\|\phi\|_1+\|\phi_x\|_1
\leq \frac{1}{\sqrt{L}}\,\|\phi\|_2+\sqrt{L}\,\|\phi_x\|_2,
\]
where the last inequality follows from the monotonicity of norms.
\item Let $\phi:\Omega\subset\mathbb{R}\rightarrow\mathbb{R}$ belong
to the class $\Sobo^{2,1}(\Omega)$.  Then
\begin{equation}
\|\phi_{x}\|_2^2\leq L\|\phi_{xx}\|_1^2 + \frac{4}{L}\|u\|_1\|\phi_{xx}\|_1.
\label{eq:result2}
\end{equation}
\emph{Proof:}  We have the identity $\int \phi_x^2{dx}\,=-\int \phi\phi_{xx}$,
true
for any function $\phi$ with periodic boundary conditions.  Using H\"older's
inequality, this becomes
\[
\|\phi_x\|_2^2\leq \|\phi\|_\infty\|\phi_{xx}\|_1.
\]
Using the relation~\eqref{eq:result1}, this becomes
\begin{align*}
\|\phi_x\|_2^2&\leq\Bigl(\frac{1}{L}\|\phi\|_1+\|\phi_x\|_1\Bigr)\|\phi_{xx}\|_1,\\
&\leq\Bigl(\frac{1}{L}\|\phi\|_1+\sqrt{L}\,\|\phi_x\|_2\Bigr)\|\phi_{xx}\|_1,
\end{align*}
which is a quadratic inequality in $\|\phi_x\|_2^2$, with solution
\[
\|\phi_x\|_2\leq\tfrac{1}{2}\Bigl(\sqrt{L}\,\|\phi_{xx}\|_1+\sqrt{L\|\phi_{xx}\|_1^2+4L^{-1}\|\phi\|_1\|\phi_{xx}\|_1}\Bigr).
\]
By sacrificing the sharpness of the bound, we obtain a simpler one,
\[
\|\phi_x\|_2^2\leq {L}\|\phi_{xx}\|_1^2 + \frac{4}{L}\|\phi\|_1\|\phi_{xx}\|_1,
\]
as required.
\end{enumerate}

\section{}
\label{apx:B}

\noindent In this appendix, we fill in the details missing from the discussion
of the regularity of solutions in Sec.~\ref{sec:regularity}
and prove the result
\[
h,c\in L^2\left(0,T_0;\Sobo^{3,2}(\Omega)\right).
\]
To begin, we notice that Sec.~\ref{sec:a_priori_bds} gives rise to
the result
\[
\|h_x\|_2,\ \|c_x\|_2 < \infty,
\]
with time-independent bounds.  Thus, using Poincar\'e's inequality, it follows
that $h,c\in \Sobo^{1,2}(\Omega)$ and, moreover, 
\[
\sup_{\left[0,T_0\right]}\|h_x\|_2,\ \sup_{\left[0,T_0\right]}\|c_x\|_2 <
\infty.
\]
From Sec.~\ref{sec:equicontinuity},
it follows that $J$ and $\mu_x$ belong to the regularity class $L^2\left(0,T_0;L^2(\Omega)\right)$,
and hence $J$, $\mu_x\in L^2\left(0,T_0;L^1(\Omega)\right)$. 
The functions $J$, $\mu$, and $\mu_x$ take the form
\[
J=h^3 h_{xxx}-h_xh^{-1}-h^2\left(h_x c_x^2+2h c_x c_{xx}\right),
\]
\[
\mu = c^3-c - h^{-1}\left(h_x c_x + hc_{xx}\right),
\]
and
\begin{equation}
\mu_x= \left(3c^2-1\right)c_x + h^{-2}h_x^2 c_x - h^{-1}h_{xx}c_x - h^{-1}h_xc_{xx}
-c_{xxx},
\label{eq:mux}
\end{equation}
respectively.
We make use of the following observations,
%The following results will help us in our demonstration,
%
%
\begin{itemize}
\item The function $h(x,t)$ is bounded from above and below,
\[
0< h_{\mathrm{min}}\leq h(x,t)\leq h_{\mathrm{max}}<\infty,
\]
and the boundedness of $c(x,t)$, $\|c\|_\infty <\infty$. 
\item Since $\mu_x\in L^2\left(0,T_0;L^2(\Omega)\right)$, it follows
that $\mu\in L^2\left(0,T_0; L^2(\Omega)\right)$, by Poincar\'e's
inequality.
\item From this it follows that $c_x h_x + hc_{xx}$ is in the class $L^2\left(0,T_0;L^2(\Omega)\right)$.
\item  Given the inequality 
\[
\|h\mu c_x\|_2^2\leq h_{\mathrm{max}}\|\mu\|_\infty^2\|c_x\|_2^2\leq
h_{\mathrm{max}}\|c_x\|_2^2\left[\frac{1}{\sqrt{L}}\,\|\mu\|_2+\sqrt{L}\,\|\mu_x\|_2\right]^2,
\]
we have the result $h_x c_x^2 + h c_x c_{xx}\in L^2\left(0,T_0;L^2(\Omega)\right)$.
\item Similarly, since $\int_0^{T_0}\|h\mu h_x\|_2^2dt<\infty$, we have the
bound
$h_x^2c_x +hh_x c_{xx}\in  L^2\left(0,T_0;L^2(\Omega)\right)$.
\end{itemize}

Now inspection of $\mu$ shows that $c_{xx}$ is in $L^2\left(0,T_0;L^1(\Omega)\right)$,
from which follows the result $h_{xxx},c_{xxx}\in L^1\left(0,T_0;L^1(\Omega)\right)$.
 By repeating the same argument, we find that $c_{xx}\in L^2\left(0,T_0;L^2(\Omega)\right)$.
  We also have the result that $\|h_x^2c_x\|_2$ is almost always bounded.
   To show that $h_{xx}\in  L^2\left(0,T_0;L^2(\Omega)\right)$,
   we take the evolution equation for $h(x,t)$, multiply it by
   $h$, and integrate, obtaining
\[
-\tfrac{1}{2}\frac{d}{dt}\int_\Omega{dx}\,h^2-3\int_\Omega{dx}\,h^2h_x^2h_{xx}-\int_\Omega{dx}\,h_xJ_0=\int_\Omega{dx}\,h^3h_{xx}^2,
\]
where $J_0=h_xh^{-1}+2h^3c_xc_{xx}+h^2h_xc_x^2$.   The time-integral
of the first term on the
left-hand side of this equation is obviously bounded in time.  Let us examine
the time-integral of the second term,
\begin{align*}
\int_0^{T_0}{dt}\int_\Omega{dx}\,h^2h_x^2h_{xx}&\leq h_{\mathrm{max}}^2\sup_{\left[0,T_0\right]}\|h_x\|_2^2\int_0^{T_0}{dt\,}\|h_{xx}\|_\infty\\
&\leq h_{\mathrm{max}}^2\sup_{\left[0,T_0\right]}\|h_x\|_2^2\int_0^{T_0}{dt}\,\left(L^{-1}\|h_{xx}\|_1+\|h_{xxx}\|_1\right)<\infty.
\end{align*}
The third term on the left-hand side is dispatched in a similar
way, so that $\int_0^{T_0}{dt}\,\|h_{xx}\|_2^2<\infty$.  We have now
shown that~$h,c\in L^2\left(0,T_0;\Sobo^{2,2}(\Omega)\right)$.  Using
this result, together with the previous facts gathered together in
this appendix, it is readily shown that~$h,c\in
L^1\left(0,T_0;\Sobo^{3,1}(\Omega)\right)$, and it follows that
\[
h,c\in L^2\left(0,T_0;\Sobo^{3,2}(\Omega)\right).
\]
To see this more clearly, we show that~$c_{xxx}$ is in the above
class.  For example, consider a typical term in~$\mu_x$,
\begin{align*}
\int_0^{T_0}{dt}\,\int_\Omega{dx}\,h_x^2c_{xx}^2%&\leq \sup_{\left[0,T_0\right]}\|h_x\|_2^2\int_0^T{dt}\,\|c_{xx}\|_\infty^2\\
&\leq \sup_{\left[0,T_0\right]}\|h_x\|_2^2\int_0^{T_0}{dt}\,\left(L^{-1}\|c_{xx}\|_1+\|c_{xxx}\|_1\right)^2<\infty,
\end{align*}
which implies~$h_x c_{xx} \in L^2\left(0,T_0;L^2(\Omega)\right)$.
This bound, together with $h\mu h_x\in
L^2\left(0,T_0;L^2(\Omega)\right)$, implies $h_x^2 c_x \in
L^2\left(0,T_0;L^2(\Omega)\right)$.  Gathering all these results, we have
\[
\mu_x\,,\ h_xc_{xx}\,,\ h_x^2c_x\,,\ h_{xx}c_x \in L^2\left(0,T_0;L^2(\Omega)\right)
\]
From~\eqref{eq:mux}, it follows that $c_{xxx}$ is in this class as
well.  A similar argument holds for $h_{xxx}$.  Thus, the solution
$\left(h,c\right)$ belongs to the regularity class
\begin{equation*}
\left(h,c\right)\in L^{\infty}\left(0,T_0;\Sobo^{1,2}(\Omega)\right)
\cap L^2\left(0,T_0;\Sobo^{3,2}(\Omega)\right)
\cap C^{\frac{1}{2},\frac{1}{8}}\left(\Omega\times\left[0,T_0\right]\right).
\end{equation*}

\section{}
\label{apx:C}
\noindent
In this appendix, we describe in full the proof of the
the uniqueness of solutions sketched in Sec.~\ref{sec:uniqueness}.
We consider two solution pairs $\left(h,c\right)$ and $\left(h',c'\right)$
and form the difference $\left(\delta h,\delta c\right)=\left(h-h',c-c'\right)$.
 Given the initial conditions $\left(\delta c(x,0),\delta h(x,0)\right)=\left(0,0\right)$,
 we show that $\left(\delta h,\delta c \right)=\left(0,0\right)$ for all
 time, that is, that the solution we have constructed is unique.
We observe that that the equation for the difference $\delta c$ can be written
in the form
\begin{equation}
\frac{\partial}{\partial t}\delta{c}+\frac{\partial^4}{\partial{x}^4}\delta{c}=\delta
\varphi\left(x,t\right),
\label{eq:delta_c_unique_appendix}
\end{equation}
where 
$\delta \varphi\left(x,t\right) \in L^2\left(0,T_0;L^2(\Omega)\right)$,
and where $\delta \varphi\left(\delta c=0\right)=0$.
As discussed previously, this equation has a unique solution, given smooth
initial data.
Since $\delta{c}=0$ satisfies Eq.~\eqref{eq:delta_c_unique_appendix}, and
since $\delta{c}\left(x,0\right)=0$, it follows that $\delta{c}=0$ for all
times $t\in\left[0,T_0\right]$.

It is now possible to formulate an equation for the difference $\delta{h}$
by subtracting the evolution equations of $h$ and $h'$ from one another,
mindful that $\delta c=0$.  We multiply the resulting equation by $\delta{h}_{xx}$
and integrate over space to obtain
\begin{multline*}
\tfrac{1}{2}\frac{d}{dt}\int_\Omega{dx}\,\delta h_x^2 +\int_\Omega{dx}\,h^3\delta
h_{xxx}^2=
\int_\Omega{dx}\,\delta h_{xxx}^2\delta h_x\left(h^{-1}+h'^2c_x^2\right)\\
-\int_\Omega{dx}\,\left(h^3-h'^3\right)\left(h'_{xxx}+4c_xc_{xx}\right)\delta
h_{xxx}\\
+ \int_\Omega{dx}\,\delta h_{xxx}\left[h'_x\left(h^{-1}-h'^{-1}\right)+h_xc_x^2\left(h^2-h'^2\right)\right].
\end{multline*}
Using the lower bound on $h(x,t)\geq h_{\mathrm{min}}>0$ and Young's
first inequality, this equation is transformed into an inequality,
\begin{multline*}
\tfrac{1}{2}\frac{d}{dt}\int_\Omega{dx}\,\delta h_x^2 +h_{\mathrm{min}}^3\int_\Omega{dx}\,\delta
h_{xxx}^2\leq
\kappa_1\int_\Omega{dx}\,\delta h_{xxx}^2+\frac{1}{4\kappa_1}\int_\Omega{dx}\,\delta{h}_x^2\left(h^{-1}+h'^2c_x^2\right)^2\\
+\kappa_2\int_\Omega{dx}\,\delta{h}_{xxx}^2+\frac{1}{4\kappa_2}\int_\Omega{dx}\,\left(h^3-h'^3\right)^2\left(h'_{xxx}+4c_xc_{xx}\right)^2\\
+
\kappa_3\int_\Omega{dx}\,\delta h_{xxx}^2+\frac{1}{4\kappa_3}\int_\Omega{dx}\,\left[h'_x\left(h^{-1}-h'^{-1}\right)+h_xc_x^2\left(h^2-h'^2\right)\right]^2,
\end{multline*}
where $\kappa_1$, $\kappa_2$, and $\kappa_3$ are arbitrary positive constants.
 By choosing $\kappa_1+\kappa_2+\kappa_3=h_{\mathrm{min}}^3$, the inequality
 simplifies to
\begin{multline*}
2\kappa\frac{d}{dt}\int_\Omega{dx}\,\delta h_x^2\leq
\int_\Omega{dx}\,\delta{h}_x^2\left(h^{-1}+h'^2c_x^2\right)^2\\+
\int_\Omega{dx}\,\delta{h}^2\left[\left(h'_{xxx}+4c_xc_{xx}\right)^2+\left(h_x'+h_xc_x^2\right)^2\right],
\end{multline*}
where $\kappa$ is another positive constant.  We integrate over the time
interval $\left[0,T\right]$, and use the fact that $\|\delta{h}_x\|_2\left(0\right)=0$
to obtain
\begin{multline*}
2\kappa\sup_{\tau\in\left[0,T\right]}\|\delta{h}_x\|_2^2\left(\tau\right)
\leq\sup_{\tau\in\left[0,T\right]}\|\delta{h}_x\|_2^2\left(\tau\right)\int_0^T{dt}\,\|h^{-1}+h'^2c_x^2\|_\infty^2\\
+\sup_{\tau\in\left[0,T\right]}\|\delta{h}\|_\infty^2\left(\tau\right)\int_0^T{dt}\int_\Omega{dx}\,\left[\left(h'_{xxx}+4c_xc_{xx}\right)^2+\left(h_x'+h_xc_x^2\right)^2\right].
\end{multline*}
The Poincar\'e inequality can be combined with the one-dimensional differential
inequalities discussed in Appendix A to yield the relation $\|f\|_\infty\leq
\kappa_P\|f_x\|_2$, where $f$ is some mean-zero function and $\kappa_P$ is
an $f$-independent constant.  We therefore arrive at
\begin{multline*}
2\kappa\sup_{\tau\in\left[0,T\right]}\|\delta{h}_x\|_2^2\left(\tau\right)
\leq\sup_{\tau\in\left[0,T\right]}\|\delta{h}_x\|_2^2\left(\tau\right)\int_0^T{dt}\,\|h^{-1}+h'^2c_x^2\|_\infty^2\\
+\kappa_P^2\sup_{\tau\in\left[0,T\right]}\|\delta{h}_x\|_2^2\left(\tau\right)\int_0^T{dt}\left(\|h'_{xxx}+4c_xc_{xx}\|_2^2+\|h_x'+h_xc_x^2\|_2^2\right).
\end{multline*}
Using the results of Sec.~\ref{sec:regularity}, it is
readily shown that $h^{-1}+h'^2c_x^2\in L^2\left(0,T;L^{\infty}(\Omega)\right)$,
and that the functions $h'_{xxx}+4c_xc_{xx}$ and $h_x'+h_xc_x^2$ belong to
the class $L^2\left(0,T;L^2(\Omega)\right)$.  By choosing $T$
sufficiently small, it is possible to impose
\[
\frac{1}{2\kappa}\left[\int_0^T{dt}\,\|h^{-1}+h'^2c_x^2\|_\infty^2+\kappa_P^2\int_0^T{dt}\left(\|h'_{xxx}+4c_xc_{xx}\|_2^2+\|h_x'+h_xc_x^2\|_2^2\right)\right]<1,
\]
which in turn forces $\sup_{\tau\in\left[0,T\right]}\|\delta{h}_x\|_2^2=0$,
and hence the solution is unique.

%%\bibliographystyle{unsrt}
%\bibliographystyle{plain}
%\bibliography{thin_film_bibliography}

\end{document}